\documentstyle{article}

\def\mj{{\mathbf{1}}}
\def\pl{\!+\!}

\def\cirk{\,{\raisebox{.3ex}{\tiny $\circ$}}\,}
\def\prop#1#2{\vspace{2ex} \noindent{\sc #1.} {\it #2} \par \vspace{2ex}}
\def\dkz{\noindent{\sc Proof. }}
\def\qed{\hfill $\dashv$}
\def\str{\rightarrow}
\def\strt{\stackrel{\textbf{.}\,}{\rightarrow}}

\begin{document}

\title{Ordinals in Frobenius Monads}
\author{\small {\sc Kosta Do\v sen} and {\sc Zoran Petri\' c}
\\[1ex]
{\small Mathematical Institute, SANU}\\[-.5ex]
{\small Knez Mihailova 36, p.f.\ 367, 11001 Belgrade,
Serbia}\\[-.5ex]
{\small email: \{kosta, zpetric\}@mi.sanu.ac.rs}}
\date{}
\maketitle

\begin{abstract}
\noindent This paper provides geometrical descriptions of the
Frobenius monad freely generated by a single object. These
descriptions are related to results connecting Frobenius algebras
and topological quantum field theories. In these descriptions,
which are based on coherence results for self-adjunctions
(adjunctions where an endofunctor is adjoint to itself), ordinals
in $\varepsilon_0$ play a prominent role. The paper ends by
considering how the notion of Frobenius algebra induces the
collapse of the hierarchy of ordinals in $\varepsilon_0$, and by
raising the question of the exact categorial abstraction of the
notion of Frobenius algebra.
\end{abstract}

\vspace{.3cm}


\vspace{.5ex}


\section{Introduction}

The purpose of this paper is to connect two seemingly distant and
unrelated topics: Frobenius algebras and ordinals contained in the
infinite denumerable ordinal $\varepsilon_0$ (namely, the least
ordinal $\xi$ such that ${\omega^\xi=\xi}$). Frobenius algebras
play an important role in topology, mathematical physics and
algebra (see \cite{K03} and references therein), while
$\varepsilon_0$ is usually deemed interesting only for
set-theorists and proof-theorists.

The categorial abstraction of the notion of Frobenius algebra
leads to the notion of Frobenius monad (for some more details, see
below). The structure of a Frobenius monad is given by a category
with an endofunctor that bears both the structure of a monad (or
triple) and a comonad, and satisfies moreover additional
conditions called Frobenius equations (see the next section).

The notion of Frobenius monad is closely related to a special kind
of adjoint situation where two functors (not necessarily distinct)
are both left and right adjoint to each other (see \cite{M65},
\cite{CMS97}, \cite{K02}, \cite{M03}, \cite{S04}, \cite{L06},
\cite{CH09}, and further references in these papers). Adjunction
is a central notion in category theory, in logic, and perhaps in
mathematics in general (see \cite{ML98} and \cite{LAW69a}), and
the connection of this notion with the notion of Frobenius monad
may serve to explain the importance of the latter.

One of the goals of this paper is to show that the notion of
adjunction where two functors are both left and right adjoint to
each other amounts, in a sense that we will make precise, to the
notion of self-adjunction, which we have investigated in
\cite{DP03}. A self-adjunction is an adjoint situation where an
endofunctor is both left and right adjoint to itself. So we find a
close relationship between Frobenius monads and self-adjunctions.
Through this relationship, we can prove coherence results for
Frobenius monads, by relying on a coherence result that we have
previously established for self-adjunctions. (That self-adjunction
arises in the context of Frobenius monads was noted in
\cite{RSW05}, Note after Definition 2.9; this is however implicit
already in \cite{LAW69}, pp.\ 151-152, in \cite{CW87}, Theorem
2.4, and in \cite{K03}, Chapter 2.)

These coherence results assert that there is a faithful functor
from the category of a freely generated Frobenius monad to
categories that serve as manageable models, which we will consider
in Section~6. This faithful functor is here an isomorphism. With
our categories that serve as models we can easily decide whether a
diagram of arrows commutes. In logical terms, this is like proving
completeness with respect to a manageable model, which helps us to
solve the decision problem. Coherence here is analogous to the
isomorphism that exists between the syntactically constructed
freely generated monad and the simplicial category (see
\cite{DP08b}, Section~3, \cite{D08}, Section~4, and references
therein).

The coherence we establish is also the gist of the connection
between the notions of Frobenius monad and two-dimensional
topological quantum field theory (2TQFT). A 2TQFT may be
understood as a functor from the category \emph{2Cob}, whose
arrows are cobordisms in dimension 2, to the category $\mbox{\it
Vect}_K$ of finite-dimensional vector spaces over the field $K$.
In terms of category theory, a Frobenius algebra is characterized
by a monoidal functor from the Frobenius monad freely generated by
a single object to $\mbox{\it Vect}_K$, modulo the strictification
of $\mbox{\it Vect}_K$ with respect to its monoidal structure
given by the tensor product and $K$ (cf.\ the beginning of
Section~7). A Frobenius algebra is the image of the object $1$ of
the Frobenius monad. The main result here is that 2TQFTs
correspond bijectively, modulo a skeletization of \emph{2Cob}, to
commutative Frobenius algebras. This result is stated officially
as a result about equivalence of categories (see \cite{K03},
Section 3.3).

An alternative result with the same mathematical content is that
the free commutative Frobenius monad is isomorphic to the skeleton
of \emph{2Cob}. From that alternative result, the former result
follows immediately. This alternative result may be conceived as a
coherence result for commutative Frobenius monads.

Our coherence results for Frobenius monads mentioned above are
more general. They deal with Frobenius monads in general, and not
only commutative ones. Because of that, infinite ordinals
contained in $\varepsilon_0$ enter into the picture. They arise
naturally in our principal category that serves as a model, which
bears some similarity to \emph{2Cob}. It is a kind of planar
version of \emph{2Cob}. Something related to this category has
been described topologically in a 2-categorial context in
\cite{KL01} (Appendix~C; see also \cite{L08}). The infinite
ordinal structure of this category is however mentioned neither in
this book, nor in the papers mentioned in the third paragraph, nor
in \cite{K03}. In \cite{K03} (Section 3.6.20) we find only the
vague conclusion that this ordinal structure, with which we want
to deal, is ``nearly about any possible drawing you can imagine''.
This structure is the main novelty we obtain when we reject
commutativity and pass to Frobenius monads in general.

This structure could be described by other means than by the
ordinals in $\varepsilon_0$. What we need is the commutative
monoid with one unary operation freely generated by the empty set
of generators (see Section~6). This monoid can be isomorphically
represented in the positive integers too, but we believe its
isomorphic representation in $\varepsilon_0$, which is quite
natural, is worth investigating.

Towards the end of his book \cite{K03} (Sections 3.6.16-27), J.
Kock discusses heuristically a project to describe geometrically
the freely generated Frobenius monad, and leaves the matter as a
challenge to the reader (Section 3.6.26). In this paper, one can
find an answer to this challenge.

To make the hierarchy of ordinals in $\varepsilon_0$ collapse, and
pass to something that amounts to \emph{2Cob}, we need not assume
commutativity. In the last two sections of this paper, we show how
the notion of Frobenius algebra requires that the notion of
Frobenius monad be extended with further assumptions, which
produce the collapse of the hierarchy. The culprit for this
collapse is the symmetry of $\mbox{\it Vect}_K$, without assuming
that the Frobenius algebra is commutative (the Frobenius objects
in symmetric monoidal categories of \cite{H04}, Section~2, involve
such a collapse too). We know that such a collapse must take
place, but we do not know what should be its exact extent. In that
context, we consider the collapse brought about by the assumption
of separability in Frobenius algebras, for which the exact
categorial abstraction is the notion of separable matrix Frobenius
monad in the last section of the paper. We leave however as an
open question what is the exact categorial abstraction of the
notion of Frobenius algebra.

This paper is a companion to \cite{DP08b}, but, except for some
side comments, an acquaintance with that paper is not
indispensable. We rely however, as we said above, on the results
of \cite{DP03}, and we assume an acquaintance with parts of that
earlier paper, though some of the essential matters we need are
reviewed in Section~6. We assume also the reader is acquainted
with some basic notions of category theory, which may all be found
in \cite{ML98}, but, for the sake of notation, we define some of
these basic notions below.

\section{The free Frobenius monad}

A \emph{Frobenius monad} is a structure made of a category $\cal
A$, an endofunctor $M$ of $\cal  A$ (i.e.\ a functor from $\cal A$
to $\cal  A$) and the natural transformations
\begin{tabbing}
\hspace{9em}\=$\varepsilon^\Box\,$\=$:M\strt I_{\cal
A}$,\hspace{7em}\= $\varepsilon^\Diamond\,$\=$:I_{\cal A}\strt M$,
\\*[1ex]
\>$\delta^\Box$\>$:M\strt MM$,\> $\delta^\Diamond$\>$:MM\strt M$,
\end{tabbing}
for $I_{\cal A}$ being the identity functor of $\cal  A$, such
that $\langle{\cal
A},M,\varepsilon^\Diamond,\delta^\Diamond\rangle$ is a monad,
$\langle{\cal A},M,\varepsilon^\Box,\delta^\Box\rangle$ is a
comonad, and, moreover, for every object $A$ of $\cal  A$, the
following \emph{Frobenius} equations hold:
\[
M\delta^\Diamond_A\cirk\delta^\Box_{MA}=\delta^\Diamond_{MA}\cirk
M\delta^\Box_A=\delta^\Box_A\cirk\delta^\Diamond_A.
\]
(For easier comparison, we use here, with slight modifications,
the notation with the modal superscripts $\Box$ and $\Diamond$,
which was introduced in \cite{DP08b}.)

The equations defining the notions of monad and comonad are given
below. For the Frobenius equations the reader may consult
\cite{K03} (in particular, Lemma 2.3.19, and \cite{DP08b},
Sections 6-7; as far as we know, and according to \cite{K06}, the
first appearance of these equations is in \cite{CW87}). Lawvere
introduced in \cite{LAW69} (pp.\ 151-152) the notion of Frobenius
monad with the equations
\[
M\varepsilon^\Box_A\cirk M\delta^\Diamond_A\cirk\delta^\Box_{MA}=
\varepsilon^\Box_{MA}\cirk\delta^\Diamond_{MA}\cirk
M\delta^\Box_A=\delta^\Diamond_A,
\]
or, alternatively, the dual equations
\[
\delta^\Diamond_{MA}\cirk M\delta^\Box_A\cirk
M\varepsilon^\Diamond_A=
M\delta^\Diamond_A\cirk\delta^\Box_{MA}\cirk\varepsilon^\Diamond_{MA}=
\delta^\Box_A,
\]
which can replace the Frobenius equations. (In the terminology of
\cite{DP08b}, Section~8, a Frobenius monad is a dyad, or codyad,
where $\Box$ and $\Diamond$ coincide.)

The category \emph{Frob} of the Frobenius monad freely generated
by a single object, denoted by $0$, has as objects the natural
numbers ${n\geq 0}$, where $n$ stands for a sequence of $n$
occurrences of $M$; so $Mn$ is ${n\pl 1}$. The arrows of this
category are defined syntactically as equivalence classes of
\emph{arrow terms}, which are defined inductively as follows. The
primitive arrow terms of \emph{Frob} are
\begin{tabbing}
\hspace{16.3em}$\mj_n\!:n\str n$,
\\[1ex]
\hspace{9em}\=$\varepsilon^\Box_n\,$\=$:n\pl 1\str
n$,\hspace{7em}\= $\varepsilon^\Diamond_n\,$\=$:n\str n\pl 1$,
\\*[1ex]
\>$\delta^\Box_n$\>$:n\pl 1\str n\pl 2$,\>
$\delta^\Diamond_n$\>$:n\pl 2\str n\pl 1$.
\end{tabbing}
The remaining arrow terms of \emph{Frob} are defined inductively
out of these with the clauses:
\begin{tabbing}
\mbox{\hspace{1.7em}}\=if ${f\!:n\str m}$ and ${g\!:m\str k}$ are
arrow terms, then so is ${(g\cirk f)\!:n\str k}$;
\\*[1ex]
\>if ${f\!:n\str m}$ is an arrow term, then so is ${Mf\!:n\pl
1\str m\pl 1}$.
\end{tabbing}
We take for granted the outermost parentheses of arrow terms, and
omit them. (Further omissions of parentheses will be permitted by
the associativity of $\cirk$.)

The least equivalence relation, congruent with respect to $\cirk$
and $M$, by which we obtain the arrows of \emph{Frob} is such
that, first, we have the \emph{categorial} equations of
composition with $\mj$ and associativity of composition $\cirk$,
and the \emph{functorial} equations for $M$ (see \cite{DP08b},
Section~2). We have next the \emph{naturality} equations:
\begin{tabbing}
\hspace{1.7em}\=($\varepsilon^\Box$~{\it
nat})\hspace{2em}\=$f\cirk\varepsilon^\Box_n=\varepsilon^\Box_m\cirk
Mf$,\hspace{4.2em}\=($\varepsilon^\Diamond$~{\it
nat})\hspace{2em}\=$\varepsilon^\Diamond_m\cirk
f=Mf\cirk\varepsilon^\Diamond_n$,
\\*[1ex]
\>($\delta^\Box$~{\it
nat})\>$MMf\cirk\delta^\Box_n=\delta^\Box_m\cirk
Mf$,\>($\delta^\Diamond$~{\it nat})\>$\delta^\Diamond_m\cirk
MMf=Mf\cirk\delta^\Diamond_n$,
\end{tabbing}
the \emph{comonad} and \emph{monad} equations:
\begin{tabbing}
\hspace{1.7em}\=$(\delta^\Box)$\hspace{3.5em}\=$M\delta^\Box_n\!$\=$\cirk\delta^\Box_n=
\delta^\Box_{n+1}\cirk\delta^\Box_n$,\hspace{2.8em}\=
$(\delta^\Diamond)$\hspace{3.6em}\=$\delta^\Diamond_n\cirk
M\delta^\Diamond_n$\=$=\delta^\Diamond_n\cirk
\delta^\Diamond_{n+1}$,
\\[1ex]
\>$(\Box\beta)$\>
$\varepsilon^\Box_{n+1}$\>$\cirk\delta^\Box_n=\mj_{n+1}$,\>
$(\Diamond\beta)$\>$\delta^\Diamond_n\cirk
\varepsilon^\Diamond_{n+1}$\>$=\mj_{n+1}$,
\\[1ex]
\>$(\Box\eta)$\>
$M\varepsilon^\Box_n$\>$\cirk\delta^\Box_n=\mj_{n+1}$,\>
$(\Diamond\eta)$\>$\delta^\Diamond_n\cirk
M\varepsilon^\Diamond_n$\>$=\mj_{n+1}$,
\end{tabbing}
and, finally, the \emph{Frobenius} equations where $A$ is replaced
by $n$. The equations $(\delta^\Box)$ and $(\delta^\Diamond)$ are
redundant in this axiomatization (see \cite{K03}, Proposition
2.3.24, and \cite{DP08b}, Section~6; they do not seem however to
be redundant when the Frobenius equations are replaced by
Lawvere's equations).

The category \emph{Frob} together with its Frobenius monad
structure is freely generated in the following sense. It is the
image of a singleton set under the left adjoint of the forgetful
functor from the category of Frobenius monads (whose arrows are
functors preserving the Frobenius monad structure on the nose,
i.e.\ exactly) to the category \emph{Set}; this forgetful functor
assigns to a Frobenius monad the set of objects of its underlying
category. Another alternative would be to take \emph{Frob}
together with its Frobenius monad structure as the image of the
trivial one-object category under the left-adjoint of the
forgetful functor from the category of Frobenius monads to the
category \emph{Cat} of small categories (whose arrows are
functors); this forgetful functor assigns to a Frobenius monad its
underlying category. As with other freely generated structures,
the syntactic construction above succeeds because the notion of
Frobenius monad is equationally presented.

The category \emph{Frob} has a strict monoidal structure. The
$\otimes$ of this monoidal structure is addition on objects. We
define ${\mj_n\otimes f}$ as $M^n f$, where $M^n$ is a sequence of
${n\geq 0}$ occurrences of $M$, while ${f\otimes\mj_n}$ is defined
by increasing the subscripts of $f$ by the natural number $n$.
Then for ${f_1\!:n_1\str m_1}$ and ${f_2\!:n_2\str m_2}$ we have
\[
f_1\otimes f_2=_{df}(f_1\otimes\mj_{m_2})\cirk(\mj_{n_1}\otimes
f_2).
\]
The category \emph{Frob} was envisaged as a monoidal category in
\cite{K03} (Section 3.6.16).

The category $\cal M$ of the monad freely generated by a single
object $0$ is defined like \emph{Frob} save that we omit the arrow
terms $\varepsilon^\Box_n$ and $\delta^\Box_n$, and whatever
involves them. By omitting $\varepsilon^\Diamond_n$ and
$\delta^\Diamond_n$, we define analogously the comonad freely
generated by $0$.

\section{Free adjunctions and monads}

An adjunction is given by two categories $\cal  A$ and $\cal  B$,
a functor $F$ from $\cal  B$ to $\cal  A$, the \emph{left
adjoint}, a functor $G$ from $\cal  A$ to $\cal  B$, the
\emph{right adjoint}, a natural transformation ${\gamma\!:I_{\cal
B}\strt GF}$, the \emph{unit} of the adjunction, and a natural
transformation ${\varphi\!:FG\strt I_{\cal A}}$, the \emph{counit}
of the adjunction, which satisfy the following \emph{triangular}
equations for every object $B$ of $\cal  B$ and every object $A$
of $\cal  A$:
\[
\varphi_{FB}\cirk
F\gamma_B=\mj_{FB},\hspace{5em}G\varphi_A\cirk\gamma_{GA}=\mj_{GA}.
\]

The adjunction freely generated by a single object $0$ on the
$\cal  B$ side is defined in syntactical terms analogously to
\emph{Frob} (see \cite{D99}, Chapter~4, for a detailed
exposition). In this free adjunction, the objects of $\cal  B$ are
$0$, $GF0$, $GFGF0$, etc., while those of $\cal  A$ are $F0$,
$FGF0$, $FGFGF0$, etc. This notion of freely generated adjunction
is essentially the same as a 2-categorial notion that may be found
in \cite{A74}, \cite{SS86} (cf.\ also \cite{KS74}) and \cite{L08}.
If we consider the sub-2-category of the 2-category \emph{Cat} of
small categories whose only 0-cells are $\cal  A$ and $\cal  B$,
whose 1-cells are made of the functors $F$ and $G$, and whose
2-cells are made of the natural transformations $\varphi$ and
$\gamma$, and other natural transformations built of these two
together with $F$ and $G$ (as, for example, $F\gamma$, $GF\gamma$,
etc.), we obtain a 2-category isomorphic in the 2-categorial sense
to the free category \emph{Ad} of \cite{A74} (called \emph{Adj} in
\cite{SS86}). This does not depend on the number of generators of
our free adjunction, provided it is not zero, and they may be
either on the $\cal  A$ side or on the $\cal B$ side.

The connection of our notion of free adjunction with the
2-category \emph{Ad} may also be construed as follows. In addition
to what we have above, we should consider the adjunction freely
generated by an object on the $\cal A$ side different from the
generating object on the $\cal B$ side, which altogether gives us
four categories with disjoint sets of objects and arrows. These
four categories are isomorphic respectively to the categories
\emph{Hom}$({\cal A},{\cal A})$, \emph{Hom}$({\cal A},{\cal B})$,
\emph{Hom}$({\cal B},{\cal B})$ and \emph{Hom}$({\cal B},{\cal
A})$ that may be found in the 2-categorial approach of \cite{A74}
and the other references above. Roughly speaking, one has only to
understand our freely generated objects as 1-cells, and add
0-cells, to pass to the 2-categorial approach. In
contradistinction to that approach, we restrict ourselves to
syntactically constructed free adjunctions within the category
\emph{Cat}, and we make explicit the free generators, but the
mathematical content is essentially the same. (The mathematical
content changes by moving to a new level of categorification with
the pseudoadjunctions of \cite{V92} and \cite{L00}.)

We give a new simple proof of the following result of \cite{A74}
(\emph{Corollaire} 2.8), which connects the category $\cal  M$ of
the free monad defined at the end of the preceding section with
the category $\cal  B$ of the adjunction freely generated by $0$
on the $\cal  B$ side. This result is interesting for us, because
it is at the base of a more complicated result concerning
\emph{Frob} that we establish in Section~5.

\prop{Proposition}{The categories $\cal  M$ and $\cal  B$ are
isomorphic.}

\dkz This isomorphism is proved syntactically by defining first by
induction a functor $I$ from $\cal  M$ to $\cal  B$ for which we
have
\begin{tabbing}
\hspace{1.7em}\=$I0=0$,\hspace{9em}\=$I(n\pl
1)=GFIn$,\hspace{6.3em}\=
\\[1ex]
\>$I\varepsilon^\Diamond_n=\gamma_{In}$,\>$I\delta^\Diamond_n=G\varphi_{FIn}$,\>$I\mj_n=\mj_{In}$,
\\[1ex]
\>$I(h_2\cirk h_1)=Ih_2\cirk Ih_1$,\>$IMh=GFIh$.
\end{tabbing}
We verify that $I$ is indeed a functor by induction on the length
of derivation of an equation of $\cal  M$.

Next we define by induction a functor $J$ from the category ${\cal
B}+{\cal A}$, which is the disjoint union of the categories $\cal
B$ and $\cal  A$ of the free adjunction, to the category~$\cal M$.
For $J$ we have
\begin{tabbing}
\hspace{1.7em}\=$J0=0$,\hspace{8.7em}\=$JGFB=JFB=JB\pl
1$,\hspace{2.9em}\=
\\[1ex]
\>$J\gamma_B=\varepsilon^\Diamond_{JB}$,\>$J\varphi_A=\delta^\Diamond_{JA-1}$,\>
$J\mj_C=\mj_{JC}$,
\\[1ex]
\>$J(h_2\cirk h_1)=Jh_2\cirk Jh_1$,\>$JGf=Jf$,\>$JFg=MJg$.
\end{tabbing}
To verify that $J$ is indeed a functor, which is done by induction
on the length of derivation of an equation, we had to define it
from ${\cal B}+{\cal A}$, but there is an obvious functor $J_{\cal
B}$ from $\cal  B$ to $\cal  M$ obtained by restricting~$J$.

It is straightforward to verify by induction on the complexity of
objects and arrow terms that $I$ and $J_{\cal B}$ are inverse to
each other. So the categories $\cal  M$ and $\cal  B$ are
isomorphic.\qed

\vspace{2ex}

\noindent A more involved, graphical, proof of this proposition
may be found in \cite{D08} (Sections 6-8).

If our free adjunction is generated by a single object on the
$\cal  A$ side, then we establish the isomorphism of $\cal  A$
with the category of the comonad freely generated by a single
object (see the end of Section~2).

\section{Bijunctions and self-adjunctions}

We call \emph{trijunction} a structure made of the categories
$\cal A$ and $\cal B$, the functor $U$ from $\cal A$ to $\cal B$,
and the functors $L$ and $R$ from $\cal B$ to $\cal A$, such that
$L$ is left adjoint to $U$, with the unit $\gamma^{\cal
B}\!:I_{\cal B}\strt UL$ and counit $\varphi^{\cal A}\!:LU\strt
I_{\cal A}$, and $R$ is right adjoint to $U$, with the unit
$\gamma^{\cal A}\!:I_{\cal A}\strt RU$ and counit $\varphi^{\cal
B}\!:UR\strt I_{\cal B}$. This notion plays an important role in
\cite{DP08b}.

We call \emph{bijunction} a trijunction where the functors $L$ and
$R$ are equal. We write in this context $P$ instead of $L$ and
$R$. (Various terms have been used for bijunctions and related
notions: in \cite{M65} one finds \emph{strongly adjoint} pairs of
functors, in \cite{CMS97} \emph{Frobenius functors}, in \cite{K02}
\emph{biadjunction}, which has already been introduced for
something else in \cite{S80}, in \cite{S04} one finds
\emph{autonomous category} and \emph{Frobenius pseudomonoid}, and
in \cite{L06} \emph{ambidextrous adjunction}.)

A \emph{self-adjunction} is an adjunction where the categories
$\cal A$ and $\cal B$ are equal, and the functors $F$ and $G$,
which are now endofunctors, are also equal. We write in this
context $\cal S$ for $\cal A$ and $\cal B$, and $F$ for both $F$
and $G$. So the unit and counit of a self-adjunction are
respectively ${\gamma\!:I_{\cal S}\strt FF}$ and
${\varphi\!:FF\strt I_{\cal S}}$. Every self-adjunction is a
bijunction. (Self-adjunction is not often mentioned in textbooks
of category theory---an exception is \cite{AHS}, Chapter~5,
Exercise 19G; the notion of self-adjoint functor of \cite{T99} is
a related but different notion.)

The bijunction freely generated by a single object $0$ on the
$\cal A$ side is defined in syntactical terms analogously to
\emph{Frob} in Section~2. The objects of the category $\cal A$ are
here $0$, $PU0$, $PUPU0$, etc., while those of $\cal B$ are $U0$,
$UPU0$, $UPUPU0$, etc.

We define analogously the free self-adjunction generated by a
single object $0$. An object of the category $\cal S$ of this
self-adjunction is of the form $F^n0$, where $F^n$ is a sequence
of ${n\geq 0}$ occurrences of $F$. We identify this object with
$n$, so that $Fn$ is ${n\pl 1}$. (One can find in \cite{DP03} a
more detailed construction of $\cal S$, which is there called
${\cal L}_c$.) The category $\cal S$ is the disjoint union of the
categories $\cal S_A$, whose objects are even, and $\cal S_B$,
whose objects are odd.

For $\cal C$ being one of the categories $\cal A$ and $\cal B$ of
the penultimate paragraph, and a subscript of one of the
categories $\cal S_A$ and $\cal S_B$ of the preceding paragraph,
we can prove the following.

\prop{Proposition}{The categories $\cal C$ and $\cal S_C$ are
isomorphic.}

\dkz We define first by induction the functors $H_{\cal C}$ from
$\cal S_C$ to $\cal C$, for $\alpha$ being $\varphi$ or~$\gamma$:
\begin{tabbing}
\hspace{1.7em}\=$H_{\cal A}0=0$,\hspace{2em}\=$H_{\cal A}(2n\pl
2)=PH_{\cal B}(2n\pl 1)$,\hspace{2.9em}\=$H_{\cal B}(2n\pl
1)=UH_{\cal A}2n$,
\\*[1ex]
\>$H_{\cal A}\alpha_{2n}=\alpha^{\cal A}_{H_{\cal A}2n}$,\>\>
$H_{\cal B}\alpha_{2n+1}=\alpha^{\cal B}_{H_{\cal B}(2n+1)}$,
\\[1ex]
\>$H_{\cal C}\mj_n=\mj_{H_{\cal C}n}$,\>\>$H_{\cal C}(h_2\cirk
h_1)=H_{\cal C}h_2\cirk H_{\cal C}h_1$,
\\[1ex]
\>$H_{\cal A}Fg=PH_{\cal B}g$,\>\>$H_{\cal B}Ff=UH_{\cal A}f$.
\end{tabbing}
Next we define by induction the functors $K_{\cal C}$ from $\cal
C$ to $\cal S_C$:
\begin{tabbing}
\hspace{1.7em}\=$H_{\cal A}0=0$,\hspace{2em}\=$H_{\cal A}(2n\pl
2)=PH_{\cal B}(2n\pl 1)$,\hspace{2.9em}\=$H_{\cal B}(2n\pl
1)=UH_{\cal A}2n$,\kill

\>$K_{\cal A}0=0$,\>$K_{\cal A}PB=K_{\cal B}B\pl 1$,\>$K_{\cal
B}UA=K_{\cal A}A\pl 1$,
\\*[1ex]
\>$K_{\cal C}\alpha^{\cal C}_C=\alpha_{K_{\cal C}C}$,
\\[1ex]
\>$K_{\cal C}\mj_C=\mj_{K_{\cal C}C}$,\>\>$K_{\cal C}(h_2\cirk
h_1)=K_{\cal C}h_2\cirk K_{\cal C}h_1$,
\\[1ex]
\>$K_{\cal A}Pg=FK_{\cal B}g$,\>\>$K_{\cal B}Uf=FK_{\cal A}f$.
\end{tabbing}

We verify by induction on the length of derivation of an equation
that $H_{\cal C}$ and $K_{\cal C}$ are indeed functors. Next we
verify by induction on the complexity of objects and arrow terms
that $H_{\cal C}$ and $K_{\cal C}$ are inverse to each other. So
the categories $\cal C$ and $\cal S_C$ are isomorphic.\qed

\section{Frobenius monads and self-adjunctions}

We want to prove the following result concerning the category
\emph{Frob} of the free Frobenius monad of Section~2 and the
category $\cal S_A$ of the free self-adjunction of the preceding
section.

\prop{Proposition}{The categories \it{Frob} and $\cal S_A$ are
isomorphic.}

\dkz We define first by induction a functor $I$ from \emph{Frob}
to $\cal S_A$:
\begin{tabbing}
\hspace{1.7em}\=$I0=0$,\hspace{7.9em}\=$I(n\pl
1)=GFIn$,\hspace{3.7em}\=\kill

\>$In=2n$,
\\*[1ex]
\>$I\varepsilon^\Box_n=\varphi_{2n}$,\>$I\delta^\Box_n=F\gamma_{2n+1}$,
\\[1ex]
\>$I\varepsilon^\Diamond_n=\gamma_{2n}$,\>$I\delta^\Diamond_n=
F\varphi_{2n+1}$,\>$I\mj_n=\mj_{2n}$,
\\[1ex]
\>$I(h_2\cirk h_1)=Ih_2\cirk Ih_1$,\>$IMh=FFIh$.
\end{tabbing}
Next we define by induction a functor $J$ from $\cal S$ to
\emph{Frob}:
\begin{tabbing}
\hspace{1.7em}\=$I0=0$,\hspace{7.9em}\=$I(n\pl
1)=GFIn$,\hspace{3em}\=\kill

\>$J2n=n$,\>$J(2n\pl 1)=n\pl 1$,
\\[1ex]
\>$J\varphi_{2n}=\varepsilon^\Box_n$,\>$J\gamma_{2n+1}=\delta^\Box_n$,
\\[1ex]
\>$J\gamma_{2n}=\varepsilon^\Diamond_n$,\>$J\varphi_{2n+1}=\delta^\Diamond_n$,\>
$J\mj_C\!=\mj_{JC}$,
\\[1ex]
\>$J(h_2\cirk h_1)=Jh_2\cirk Jh_1$,\>$JFg=Jg$, for $g$ in $\cal
S_B$,\>$JF\!f\!=MJf$, for $f$ in $\cal S_A$.
\end{tabbing}

We verify by induction on the length of derivation of an equation
that $I$ and $J$ are indeed functors. We will not dwell on that
verification for $I$, while for $J$ we have to verify first that
\[
J(h\cirk\varphi_n)=J(\varphi_m\cirk FFh).
\]
If $h$ is from $\cal S_A$, then we use the equation
($\varepsilon^\Box$~{\it nat}) of Section~2. If $h$ is from $\cal
S_B$, then we proceed by induction on the complexity of $h$, by
using the Frobenius equations and the equations
$(\delta^\Diamond)$ and ($\delta^\Diamond$~{\it nat}) of
Section~2. Note that if $h$ is from $\cal S_B$, then $Jh$ can be
neither $\varepsilon^\Box_k$ nor $\varepsilon^\Diamond_k$. We
proceed analogously for
\[
J(\gamma_m\cirk h)=J(FFh\cirk \gamma_n).
\]
To verify ${Jh_1=Jh_2}$ for ${h_1=h_2}$ a triangular equation, we
use the equations $(\Box\beta)$, $(\Box\eta)$, $(\Diamond\beta)$
and $(\Diamond\eta)$ of Section~2.

There is an obvious functor $J_{\cal A}$ from $\cal S_A$ to
\emph{Frob} obtained by restricting $J$, and it is straightforward
to verify by induction on the complexity of arrow terms that $I$
and $J_{\cal A}$ are inverse to each other. So the categories
\emph{Frob} and $\cal S_A$ are isomorphic.$\dashv$

\vspace{2ex}

From this proposition and from the Proposition of the preceding
section we can conclude that \emph{Frob} is isomorphic to the
category $\cal A$ of the bijunction freely generated by a single
object on the $\cal A$ side, but the isomorphism we have
established in this section is more interesting for us, as it will
become clear in the next section.

\section{Coherence for Frobenius monads}

Out of the category $\cal S$ of the free self-adjunction of
Section~4, we build a monoid ${\cal S}^\ast$ (which in
\cite{DP03}, Section on ${\cal L}_c$ \emph{and} ${\cal L}_\omega$,
is called ${\cal L}_c^\ast$, while $\cal S$ is called ${\cal
L}_c$). On the arrows of $\cal S$, we define a total binary
operation $\ast$ based on composition in the following manner: for
$f\!:m\str n$ and $g\!:k\str l$,
\[
g\ast f=_{df}\left\{
\begin{array}{ll}
g\cirk F^{k-n}f & {\mbox{\rm if }}n\leq k
\\
F^{n-k}g\cirk f & {\mbox{\rm if }}k\leq n.
\end{array}
\right.
\]

Next, let $f\equiv g$ if and only if for some $k$ and $l$ we have
$F^kf=F^lg$ in $\cal S$. It is easy to check that $\equiv$ is an
equivalence relation on the arrows of $\cal S$, which satisfies
moreover
\[
{\mbox{\rm if }}f_1\equiv f_2\; {\mbox{\rm and }} g_1\equiv g_2,
\; {\mbox{\rm then }}g_1\ast f_1\equiv g_2\ast f_2.
\]

For every arrow $f$ of $\cal S$, let $[f]$ be $\{g\mid f\equiv
g\}$, and let ${\cal S}^\ast$ be $\{[f]\mid f$ is an arrow of
$\cal S\}$. With
\[
\begin{array}{l}
\mj=_{df}[\mj_0],
\\[.1cm]
[g][f]=_{df}[g\ast f],
\end{array}
\]
we can check that ${\cal S}^\ast$ is indeed a monoid, which in
\cite{DP03} is shown to be isomorphic to the monoid ${\cal
L}_\omega$. The monoid ${\cal S}^\ast$ is built of syntactical
material, coming from the category $\cal S$, which is presented by
generators and equations. The monoid ${\cal L}_\omega$ is just
another presentation by generators and equations of ${\cal
S}^\ast$. A reason for introducing it in \cite{DP03} was to make
simpler reduction to normal form, without being encumbered by the
sources and targets of arrow terms. We presuppose the reader is
acquainted with ${\cal L}_\omega$, but indications about how this
monoid is presented will be given below when we deal with
composition in \emph{Frobse}. This monoid interests us here only
as an auxiliary, leading towards the geometric categories
\emph{Frz} and \emph{Frobse}, which we will consider below.

Out of the material of the monoid ${\cal S}^{\ast}$ we return to
$\cal S$, by building a category isomorphic to $\cal S$, in the
following way. Let ${\cal S}^{\ast t}$ be the category whose
objects are the natural numbers, and whose arrows are the triples
${\langle [f], n, m\rangle}$ such that there is an arrow
${g:\!n\str m}$ in $[f]$; the arrow ${\langle [f], n, m\rangle}$
is of \emph{type} $n\str m$ in ${\cal S}^{\ast t}$, which means
that its source is $n$ and its target $m$. The composition of
${\langle [f_1], n, m\rangle}$ and ${\langle [f_2],m,k\rangle}$ is
defined as ${\langle [f_2][f_1],n,k\rangle}$, and the identity
arrow on $n$ is ${\langle [\mj_0],n,n\rangle}$. We can prove the
following.

\prop{Proposition 1}{The categories ${\cal S}^{\ast t}$ and $\cal
S$ are isomorphic.}

The functor from $\cal S$ to ${\cal S}^{\ast t}$ giving this
isomorphism is identity on objects and an arrow $f\!:n\str m$ is
mapped to ${\langle [f], n, m\rangle}$. To prove the proposition,
we rely on the fact that for every element $[f]$ of ${\cal
S}^\ast$, and every arrow ${g:\!n\str m}$ in $[f]$, the arrow $g$
is the only arrow in $[f]$ of the type ${n\str m}$. To establish
this fact we need to establish the following.

\prop{$\cal S$ Cancellation Lemma}{In $\cal S$, if $Ff=Fg$, then
$f=g$.}

An analogous lemma (called the $\cal L$ Cancellation Lemma) was
proved in \cite{DP03}, but with a restriction on the type of $f$
and $g$. It was stated there that the restriction can be lifted,
and it was suggested how to achieve that. The suggestion envisaged
two ways, one of which is rather straightforward (both of these
ways are however lengthy, and for lack of space a detailed
exposition was omitted.) As a matter of fact, there is a direct
way to prove the $\cal S$ Cancellation Lemma along the lines of
the ${\cal K}_c$ Cancellation Lemma of \cite{DP03}. Here are
indications for this direct proof (which presuppose an
acquaintance with \cite{DP03}).

\vspace{2ex}

\noindent {\sc Proof of the $\cal S$ Cancellation Lemma.} We
proceed as in the proof of the ${\cal K}_c$ Cancellation Lemma of
\cite{DP03} until we reach the case that $f=\varphi_0\cirk f'$ and
$g=\varphi_0\cirk g'$. Then we must ensure that $\varphi_0$ in
$\varphi_0\cirk f'$ and $\varphi_0$ in $\varphi_0\cirk g'$ are
tied in $\delta(\psi(f))$, which is equal to $\delta(\psi(g))$, to
circles encompassing the same circular forms. It is always
possible to achieve that. We conclude that
$\delta(\psi(f'))\cong_{\cal L} \delta(\psi(g'))$. This is because
$\delta(\psi(\varphi_0\cirk f'))$ and $\delta(\psi(\varphi_0\cirk
g'))$ are both $\cal L$-equivalent to $\delta(c_1^\alpha)$, for
some $\alpha\geq 1$, while $\delta(\psi(f'))$ and
$\delta(\psi(g'))$ must both be $\cal L$-equivalent to the same
$\delta(b_1^{\beta}c_1^{\gamma})$, for $\beta$ and $\gamma$ lesser
than $\alpha$. We conclude that ${f'\equiv g'}$, and since $f'$
and $g'$ are of type $0\str 2$, by the $\cal L$ Cancellation Lemma
of \cite{DP03}, we have that $f'=g'$ in $\cal S$. From that we
obtain that $f=g$ in $\cal S$.\qed

\vspace{2ex}

As $\cal S$, the isomorphic category ${\cal S}^{\ast t}$ is the
disjoint union of two categories ${\cal S}_{\cal A}^{\ast t}$ and
${\cal S}_{\cal B}^{\ast t}$, isomorphic respectively to $\cal
S_A$ and $\cal S_B$. If ${\langle [f], n, m\rangle}$ is in ${\cal
S}_{\cal A}^{\ast t}$, then $n$ and $m$ are even, and if it is in
${\cal S}_{\cal B}^{\ast t}$, then they are odd. We will now
define a category ${\cal L}^{\cal A}_\omega$, isomorphic to ${\cal
S}_{\cal A}^{\ast t}$. This category, made out of the material of
the monoid ${\cal L}_\omega$, interests us only as a stepping
stone towards the isomorphic geometric categories \emph{Frz} and
\emph{Frobse}, which we will consider in a moment.

The objects of the category ${\cal L}^{\cal A}_\omega$ are again
the natural numbers. An arrow of this category will be obtained
from an arrow ${\langle [f], n, m\rangle}$ of ${\cal S}_{\cal
A}^{\ast t}$ by replacing the class $[f]$ by the corresponding
element $e$ of ${\cal L}_\omega$, and $n$ and $m$ by respectively
$n/2$ and $m/2$; the type of ${\langle e,n/2,m/2\rangle}$ in
${\cal L}^{\cal A}_\omega$ is ${n/2\str m/2}$. We divide the
numbers in the types by two to obtain types that will correspond
to the types in \emph{Frob} (see the preceding section).

In this way, with every element of ${\cal L}_\omega$ we associate
in ${\cal L}^{\cal A}_\omega$ a denumerable infinity of types. The
generator $a^\alpha_{2n+1}$ of ${\cal L}_\omega$ (see \cite{DP03},
Section on \emph{Normal forms in} ${\cal L}_\omega$) will have as
associated types $n\pl l\pl 1\str n\pl l$, for every $l\geq 0$,
while the generator $b^\alpha_{2n+1}$ will have $n\pl l\str n\pl
l\pl 1$, and the generator $c^\alpha_{2n+1}$ will have $n\pl l\str
n\pl l$. (This typing is explained by the typing of the friezes
below.) The generator $a^\alpha_{2n+2}$ will have as associated
types $n\pl l\pl 2\str n\pl l\pl 1$, the generator
$b^\alpha_{2n+2}$ will have $n\pl l\pl 1\str n\pl l\pl 2$, and the
generator $c^\alpha_{2n+2}$ will have $n\pl l\pl 1\str n\pl l\pl
1$. Multiplication of terms now becomes composition, and takes the
types into account. Two typed terms of ${\cal L}_\omega$ stand for
the same arrow of the category ${\cal L}^{\cal A}_\omega$ if and
only if they are of the same type and equal in ${\cal L}_\omega$.
We can then assert the following.

\prop{Proposition 2}{The categories \it{Frob} and ${\cal L}^{\cal
A}_\omega$ are isomorphic.}

\noindent This follows immediately from the isomorphism of
\emph{Frob} and $\cal S_A$, proved in the preceding section, the
isomorphism of $\cal S_A$ and ${\cal S}_{\cal A}^{\ast t}$, which
follows from Proposition 1, and the isomorphism of ${\cal S}_{\cal
A}^{\ast t}$ and ${\cal L}^{\cal A}_\omega$.

From \cite{DP03} (Section on ${\cal L}_\omega$, ${\cal K}_\omega$
\emph{and friezes}) one can infer that the category ${\cal
L}^{\cal A}_\omega$ is isomorphic to a category \emph{Frz} whose
arrows are diagrams called \emph{friezes} with associated types.
Roughly speaking, a frieze is a tangle without crossings in whose
regions we find circular forms that correspond bijectively to the
ordinals contained in the infinite ordinal $\varepsilon_0$. In
\cite{DP03} one can find a proof that ${\cal L}_\omega$ is
isomorphic to the monoid of friezes, and from that the isomorphism
of the categories ${\cal L}^{\cal A}_\omega$ and \emph{Frz}
follows. So, by Proposition 2, the categories \emph{Frob} and
\emph{Frz} are isomorphic. By this last isomorphism, the arrows on
the left are mapped to the friezes on the right, with the type
associated to the friezes being those of the arrows:
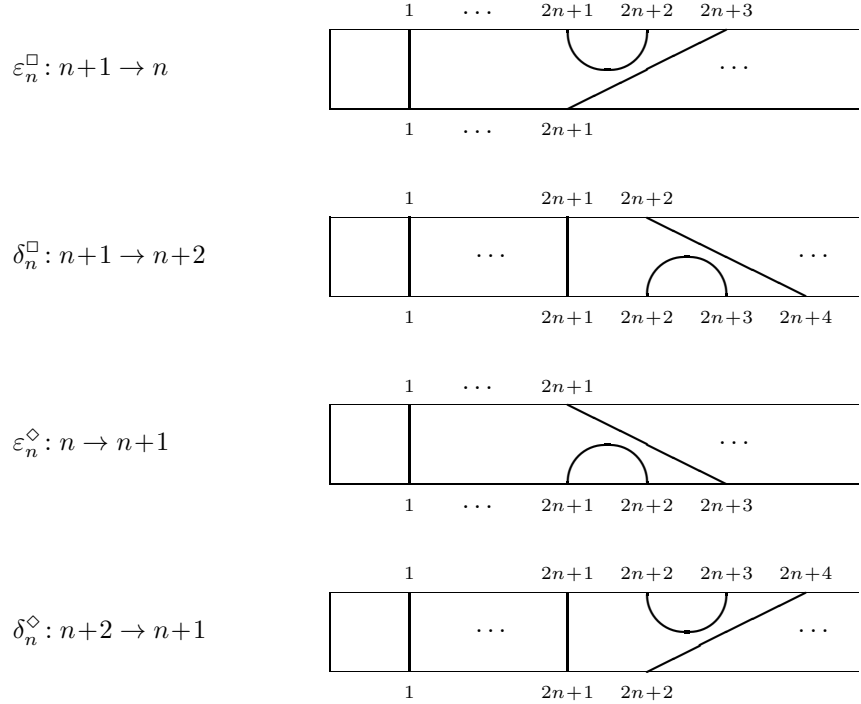
\begin{center}
\begin{picture}(200,60)(-60,0)

\put(-120,25){\makebox(0,0)[l]{${\varepsilon^\Box_n\!:n\pl 1\str
n}$}}

{\linethickness{0.05pt} \put(0,10){\line(1,0){200}}
\put(0,10){\line(0,1){30}} \put(0,40){\line(1,0){200}}}

\put(30,0){\makebox(0,0)[b]{\scriptsize $1$}}
\put(90,0){\makebox(0,0)[b]{\scriptsize $2n\pl 1$}}
\put(30,45){\makebox(0,0)[b]{\scriptsize $1$}}
\put(90,45){\makebox(0,0)[b]{\scriptsize $2n\pl 1$}}
\put(120,45){\makebox(0,0)[b]{\scriptsize $2n\pl 2$}}
\put(150,45){\makebox(0,0)[b]{\scriptsize $2n\pl 3$}}

\thicklines \put(30,10){\line(0,1){30}}
\put(90,10){\line(2,1){60}}

\put(105,40){\oval(30,30)[b]}

\put(50,1){$\ldots$} \put(50,46){$\ldots$} \put(147,25){$\ldots$}

\end{picture}
\begin{picture}(200,70)(-60,0)

\put(-120,25){\makebox(0,0)[l]{${\delta^\Box_n\!:n\pl 1\str n\pl
2}$}}

{\linethickness{0.05pt} \put(0,10){\line(1,0){200}}
\put(0,10){\line(0,1){30}} \put(0,40){\line(1,0){200}}}

\put(30,0){\makebox(0,0)[b]{\scriptsize $1$}}
\put(90,0){\makebox(0,0)[b]{\scriptsize $2n\pl 1$}}
\put(120,0){\makebox(0,0)[b]{\scriptsize $2n\pl 2$}}
\put(150,0){\makebox(0,0)[b]{\scriptsize $2n\pl 3$}}
\put(180,0){\makebox(0,0)[b]{\scriptsize $2n\pl 4$}}

\put(30,45){\makebox(0,0)[b]{\scriptsize $1$}}
\put(90,45){\makebox(0,0)[b]{\scriptsize $2n\pl 1$}}
\put(120,45){\makebox(0,0)[b]{\scriptsize $2n\pl 2$}}

\thicklines \put(30,10){\line(0,1){30}}
\put(90,10){\line(0,1){30}} \put(180,10){\line(-2,1){60}}

\put(135,10){\oval(30,30)[t]}

\put(55,25){$\ldots$} \put(177,25){$\ldots$}

\end{picture}
\begin{picture}(200,70)(-60,0)

\put(-120,25){\makebox(0,0)[l]{${\varepsilon^\Diamond_n\!:n\str
n\pl 1}$}}

{\linethickness{0.05pt} \put(0,10){\line(1,0){200}}
\put(0,10){\line(0,1){30}} \put(0,40){\line(1,0){200}}}

\put(30,45){\makebox(0,0)[b]{\scriptsize $1$}}
\put(90,45){\makebox(0,0)[b]{\scriptsize $2n\pl 1$}}
\put(30,0){\makebox(0,0)[b]{\scriptsize $1$}}
\put(90,0){\makebox(0,0)[b]{\scriptsize $2n\pl 1$}}
\put(120,0){\makebox(0,0)[b]{\scriptsize $2n\pl 2$}}
\put(150,0){\makebox(0,0)[b]{\scriptsize $2n\pl 3$}}

\thicklines \put(30,10){\line(0,1){30}}
\put(90,40){\line(2,-1){60}}

\put(105,10){\oval(30,30)[t]}

\put(50,1){$\ldots$} \put(50,46){$\ldots$} \put(147,25){$\ldots$}

\end{picture}
\begin{picture}(200,70)(-60,0)

\put(-120,25){\makebox(0,0)[l]{${\delta^\Diamond_n\!:n\pl 2\str
n\pl 1}$}}

{\linethickness{0.05pt} \put(0,10){\line(1,0){200}}
\put(0,10){\line(0,1){30}} \put(0,40){\line(1,0){200}}}

\put(30,45){\makebox(0,0)[b]{\scriptsize $1$}}
\put(90,45){\makebox(0,0)[b]{\scriptsize $2n\pl 1$}}
\put(120,45){\makebox(0,0)[b]{\scriptsize $2n\pl 2$}}
\put(150,45){\makebox(0,0)[b]{\scriptsize $2n\pl 3$}}
\put(180,45){\makebox(0,0)[b]{\scriptsize $2n\pl 4$}}

\put(30,0){\makebox(0,0)[b]{\scriptsize $1$}}
\put(90,0){\makebox(0,0)[b]{\scriptsize $2n\pl 1$}}
\put(120,0){\makebox(0,0)[b]{\scriptsize $2n\pl 2$}}

\thicklines \put(30,10){\line(0,1){30}}
\put(90,10){\line(0,1){30}} \put(180,40){\line(-2,-1){60}}

\put(135,40){\oval(30,30)[b]}

\put(55,25){$\ldots$} \put(177,25){$\ldots$}

\end{picture}
\end{center}
When ${n=0}$, the vertical thread connecting 1 at the top with 1
at the bottom does not exist in the first and the third frieze.
Note that our friezes are ``thin'' tangles that may be conceived
as the boundaries of the corresponding \emph{thick} tangles of
\cite{KL01}.

A \emph{circular form} is a finite collection of nonintersecting
circles in the plane factored through homeomorphisms of the plane
mapping one collection into another (see the definition of $\cal
L$-equivalence of friezes in \cite{DP03}, Section on
\emph{Friezes}). The circular forms obtained by composing friezes
are coded by the ordinals contained in $\varepsilon_0$ in the
following way. The circular form consisting of no circles is coded
by 0. If the circular forms $c_1$, $c_2$ and $c$ are coded by the
ordinals $\alpha_1$, $\alpha_2$ and $\alpha$ respectively, then
the circular form ${c_1c_2}$ (the disjoint union of $c_1$ and
$c_2$) is coded by the natural sum ${\alpha_1\sharp\,\alpha_2}$,
and the circular form
\begin{picture}(12,7) \put(5,2){\circle{10}}
\put(5,0){\makebox(0,0)[b]{$c$}}
\end{picture} ($c$ inside a new circle) is coded by $\omega^\alpha$. So a single circle is
coded by $\omega^0$, which is equal to 1 (see \cite{DP03}, Section
on \emph{Finite multisets, circular forms and ordinals}).

Let $\cal F$ be the commutative monoid with one unary operation
freely generated by the empty set of generators. The elements of
$\cal F$ may be identified with the hierarchy of finite multisets
obtained by starting from the empty multiset as the only
urelement, or by finite nonplanar trees with arbitrary finite
branching, or by circular forms. A monoid isomorphic to $\cal F$
is the commutative monoid
${\langle\varepsilon_0,\sharp,0,\omega^{-}\rangle}$ where $\sharp$
is binary natural sum, and we have the additional unary operation
$\omega^{-}$ (for more details on these matters, see \cite{DP03}).
Note that though the elements of $\varepsilon_0$ greater than or
equal to $\omega$ are associated with infinite ordinals, they may
be used to code finite objects, such as circular forms. Another
monoid isomorphic to $\cal F$ is the commutative monoid
${\langle\textbf{N}^+,\cdot,1,p_{\underline{\;\;}}\,\rangle}$
where $\textbf{N}^+$ is the set of natural numbers greater than 0,
the operation $\cdot$ is multiplication, and $p_n$ is the $n$-th
prime number (we are indebted for this remark to a suggestion of
Marko Sto\v si\' c).

The isomorphism of \emph{Frob} with \emph{Frz} may be understood
as a geometrical description of \emph{Frob}. Towards the end of
his book \cite{K03} (Sections 3.6.20 ff), Kock was looking for
such a description, but not exactly in the same direction. The
category \emph{Frobse}, isomorphic to \emph{Frz}, which we will
consider below, gives another alternative approach to the
geometrization of \emph{Frob} sought by Kock.

The isomorphism of \emph{Frob} and \emph{Frz} may be understood
also as a coherence result, which provides a decision procedure
for equality of arrows in \emph{Frob}. This decision procedure
involves a syntactical description of friezes given by the monoid
${\cal L}_\omega$ of \cite{DP03}, and a reduction to normal form.

Instead of the category \emph{Frz}, one can use an alternative
isomorphic category, which we will call \emph{Frobse}. In the
arrows of this category, the regions of friezes stand for
equivalence classes of an equivalence relation whose domain is
split into a source part and a target part, which are both copies
of $\textbf{N}^+$. Such equivalence relations, called \emph{split
equivalences}, are studied in \cite{DP03a} and \cite{DP09}. Split
equivalences are related to \emph{cospans} in the base category
\emph{Set} (see \cite{ML98}, XII.7, and \cite{RSW05}, Example
2.4), but unlike cospans they do not register the common target of
the two arrows making the cospan.

The split equivalences we envisage for \emph{Frobse} are
\emph{nonintersecting} in the following sense. Let the source and
target elements be identified respectively with the positive and
negative integers (so 0 does not correspond to any element). For
$a,b,c,d\in\textbf{Z}-\{0\}$, we say that ${(a,b)}$
\emph{intersects} ${(c,d)}$  when either $a<c<b<d$ or $c<a<d<b$.
An equivalence relation on ${\textbf{Z}-\{0\}}$ is
\emph{nonintersecting} when if $a$ and $b$ are in one equivalence
class, while $c$ and $d$ are in another equivalence class, then
${(a,b)}$ does not intersect ${(c,d)}$. (This is related to the
\emph{nonoverlapping} segments of \cite{DP03}, Section on
\emph{Friezes}.)

For example, instead of the frieze on the left-hand side, which is
an arrow of \emph{Frz} of the type ${2\pl l\str 1\pl l}$, we have
the nonintersecting split equivalence on the right-hand side,
which is an arrow of \emph{Frobse} of the same type:
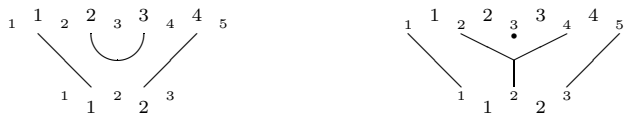
\begin{center}
\begin{picture}(200,40)

\put(0,35){\makebox(0,0)[b]{\scriptsize $1$}}
\put(20,35){\makebox(0,0)[b]{\scriptsize $2$}}
\put(40,35){\makebox(0,0)[b]{\scriptsize $3$}}
\put(60,35){\makebox(0,0)[b]{\scriptsize $4$}}

\put(20,0){\makebox(0,0)[b]{\scriptsize $1$}}
\put(40,0){\makebox(0,0)[b]{\scriptsize $2$}}

\put(-10,32){\makebox(0,0)[b]{\tiny $1$}}
\put(10,32){\makebox(0,0)[b]{\tiny $2$}}
\put(30,32){\makebox(0,0)[b]{\tiny $3$}}
\put(50,32){\makebox(0,0)[b]{\tiny $4$}}
\put(70,32){\makebox(0,0)[b]{\tiny $5$}}

\put(10,8){\makebox(0,0)[t]{\tiny $1$}}
\put(30,8){\makebox(0,0)[t]{\tiny $2$}}
\put(50,8){\makebox(0,0)[t]{\tiny $3$}}

\put(20,10){\line(-1,1){20}} \put(40,10){\line(1,1){20}}

\put(30,30){\oval(20,20)[b]}

\put(150,35){\makebox(0,0)[b]{\scriptsize $1$}}
\put(170,35){\makebox(0,0)[b]{\scriptsize $2$}}
\put(190,35){\makebox(0,0)[b]{\scriptsize $3$}}
\put(210,35){\makebox(0,0)[b]{\scriptsize $4$}}

\put(170,0){\makebox(0,0)[b]{\scriptsize $1$}}
\put(190,0){\makebox(0,0)[b]{\scriptsize $2$}}

\put(140,32){\makebox(0,0)[b]{\tiny $1$}}
\put(160,32){\makebox(0,0)[b]{\tiny $2$}}
\put(180,32){\makebox(0,0)[b]{\tiny $3$}}
\put(200,32){\makebox(0,0)[b]{\tiny $4$}}
\put(220,32){\makebox(0,0)[b]{\tiny $5$}}

\put(160,8){\makebox(0,0)[t]{\tiny $1$}}
\put(180,8){\makebox(0,0)[t]{\tiny $2$}}
\put(200,8){\makebox(0,0)[t]{\tiny $3$}}

\put(160,10){\line(-1,1){20}} \put(200,10){\line(1,1){20}}
\put(180,10){\line(0,1){10}} \put(180,20){\line(-2,1){20}}
\put(180,20){\line(2,1){20}}

\put(180,29){\circle*{2}}

\end{picture}
\end{center}
The thick white regions on the left-hand side become thin black
equivalence classes on the right-hand side, and the thin black
threads on the left-hand side become white regions on the
right-hand side. We will not obtain in this way on the right-hand
side every nonintersecting split equivalence.

The equivalence classes of those we obtain satisfy some additional
conditions. First, they are all finite, and all but finitely many
of them are such that they have just two elements---one at the top
and one at the bottom. Secondly, they are either \emph{even} or
\emph{odd}, depending on whether their members are even or odd; we
have only such even and odd equivalence classes. Finally, two
classes of the same parity cannot be immediate neighbours in the
following sense. The classes $A$ and $B$ are \emph{immediate
neighbours} when for every ${a\in A}$ and every ${b\in B}$ and
every class $C$ and every ${c_1,c_2\in C}$, if ${(a,b)}$
intersects ${(c_1,c_2)}$, then $C$ is either $A$ or $B$. The
nonintersecting split equivalences that satisfy these additional
conditions concerning their equivalence classes will be called
\emph{maximal} split equivalences.

Note that in maximal split equivalences the odd equivalence
classes are completely determined by the even equivalence classes,
and vice versa. We cannot however reject either of them because of
the ordinals. In the regions of friezes one finds finitely many
circular forms that correspond to ordinals in $\varepsilon_0$, and
we will assign these ordinals to the equivalence classes of
maximal split equivalences.

Maximal split equivalences together with a function assigning
ordinals in $\varepsilon_0$ to the equivalence classes, so that
all but finitely many have zero as value, will be called
\emph{Frobenius} split equivalences. Frobenius split equivalences
with types associated to them are the arrows of \emph{Frobse}. For
example, to the frieze on the left-hand side we assign the
Frobenius split equivalence on the right-hand side:
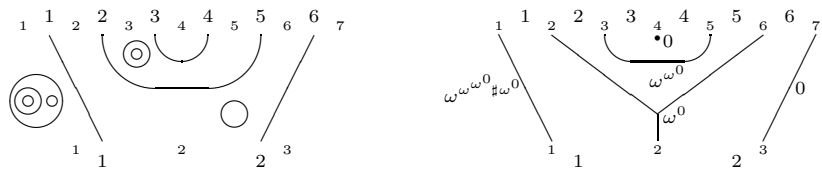
\begin{center}
\begin{picture}(280,60)

\put(0,55){\makebox(0,0)[b]{\scriptsize $1$}}
\put(20,55){\makebox(0,0)[b]{\scriptsize $2$}}
\put(40,55){\makebox(0,0)[b]{\scriptsize $3$}}
\put(60,55){\makebox(0,0)[b]{\scriptsize $4$}}
\put(80,55){\makebox(0,0)[b]{\scriptsize $5$}}
\put(100,55){\makebox(0,0)[b]{\scriptsize $6$}}

\put(20,0){\makebox(0,0)[b]{\scriptsize $1$}}
\put(80,0){\makebox(0,0)[b]{\scriptsize $2$}}

\put(-10,52){\makebox(0,0)[b]{\tiny $1$}}
\put(10,52){\makebox(0,0)[b]{\tiny $2$}}
\put(30,52){\makebox(0,0)[b]{\tiny $3$}}
\put(50,52){\makebox(0,0)[b]{\tiny $4$}}
\put(70,52){\makebox(0,0)[b]{\tiny $5$}}
\put(90,52){\makebox(0,0)[b]{\tiny $6$}}
\put(110,52){\makebox(0,0)[b]{\tiny $7$}}

\put(10,8){\makebox(0,0)[t]{\tiny $1$}}
\put(50,8){\makebox(0,0)[t]{\tiny $2$}}
\put(90,8){\makebox(0,0)[t]{\tiny $3$}}

\put(20,10){\line(-1,2){20}} \put(80,10){\line(1,2){20}}

\put(50,50){\oval(20,20)[b]} \put(50,50){\oval(60,40)[b]}

\put(33,43){\circle{10}} \put(33,43){\circle{4}}
\put(70,20){\circle{10}} \put(-5,25){\circle{20}}
\put(1,25){\circle{4}} \put(-8,25){\circle{10}}
\put(-8,25){\circle{4}}

\put(180,55){\makebox(0,0)[b]{\scriptsize $1$}}
\put(200,55){\makebox(0,0)[b]{\scriptsize $2$}}
\put(220,55){\makebox(0,0)[b]{\scriptsize $3$}}
\put(240,55){\makebox(0,0)[b]{\scriptsize $4$}}
\put(260,55){\makebox(0,0)[b]{\scriptsize $5$}}
\put(280,55){\makebox(0,0)[b]{\scriptsize $6$}}

\put(200,0){\makebox(0,0)[b]{\scriptsize $1$}}
\put(260,0){\makebox(0,0)[b]{\scriptsize $2$}}

\put(170,52){\makebox(0,0)[b]{\tiny $1$}}
\put(190,52){\makebox(0,0)[b]{\tiny $2$}}
\put(210,52){\makebox(0,0)[b]{\tiny $3$}}
\put(230,52){\makebox(0,0)[b]{\tiny $4$}}
\put(250,52){\makebox(0,0)[b]{\tiny $5$}}
\put(270,52){\makebox(0,0)[b]{\tiny $6$}}
\put(290,52){\makebox(0,0)[b]{\tiny $7$}}

\put(190,8){\makebox(0,0)[t]{\tiny $1$}}
\put(230,8){\makebox(0,0)[t]{\tiny $2$}}
\put(270,8){\makebox(0,0)[t]{\tiny $3$}}

\put(190,10){\line(-1,2){20}} \put(270,10){\line(1,2){20}}
\put(230,10){\line(0,1){10}} \put(230,20){\line(-4,3){40}}
\put(230,20){\line(4,3){40}} \put(230,50){\oval(40,20)[b]}
\put(230,49){\circle*{2}}

\put(179,30){\makebox(0,0)[r]{\scriptsize
$\omega^{\omega^{\omega^0}\!\sharp \omega^0}$}}
\put(232,23){\makebox(0,0)[tl]{\scriptsize $\omega^0$}}
\put(234,38.5){\makebox(0,0)[t]{\scriptsize $\omega^{\omega^0}$}}
\put(232,50){\makebox(0,0)[tl]{\scriptsize $0$}}
\put(282,30){\makebox(0,0)[l]{\scriptsize $0$}}

\end{picture}
\end{center}

All the  Frobenius split equivalences are generated by composition
from the following \emph{generating} Frobenius split equivalences,
which are correlated with the elements of the monoid ${\cal
L}_\omega$ mentioned on the left of the following pictures (see
\cite{DP03}, Section on \emph{Normal forms in} ${\cal L}_\omega$
\emph{and} ${\cal K}_\omega$), where we omit mentioning that an
equivalence class bears $0$; here, ${k\geq 1}$ and
${\alpha,\beta\in\varepsilon_0}$:
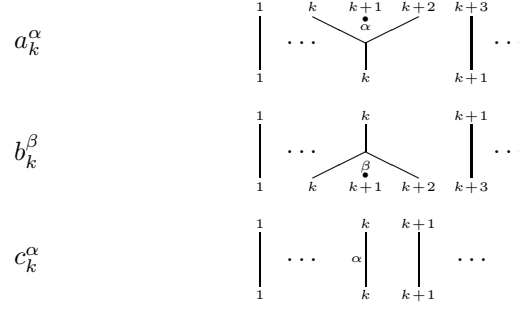
\begin{center}
\begin{picture}(310,40)

\put(22,20){\makebox(0,0)[l]{$a^\alpha_k$}}

\put(125,20){\makebox(0,0)[l]{$\ldots$}}
\put(215,20){\makebox(0,0)[r]{$\ldots$}}

\put(115,32){\makebox(0,0)[b]{\tiny $1$}}
\put(135,32){\makebox(0,0)[b]{\tiny $k$}}
\put(155,32){\makebox(0,0)[b]{\tiny $k\pl 1$}}
\put(175,32){\makebox(0,0)[b]{\tiny $k\pl 2$}}
\put(195,32){\makebox(0,0)[b]{\tiny $k\pl 3$}}

\put(115,5){\makebox(0,0)[b]{\tiny $1$}}
\put(155,5){\makebox(0,0)[b]{\tiny $k$}}
\put(195,5){\makebox(0,0)[b]{\tiny $k\pl 1$}}

\put(115,10){\line(0,1){20}} \put(155,10){\line(0,1){10}}
\put(155,20){\line(-2,1){20}} \put(155,20){\line(2,1){20}}
\put(195,10){\line(0,1){20}}

\put(155,27){\makebox(0,0)[t]{\tiny $\alpha$}}

\put(155,29){\circle*{2}}

\end{picture}

\begin{picture}(310,40)

\put(22,20){\makebox(0,0)[l]{$b^\beta_k$}}

\put(125,20){\makebox(0,0)[l]{$\ldots$}}
\put(215,20){\makebox(0,0)[r]{$\ldots$}}

\put(115,5){\makebox(0,0)[b]{\tiny $1$}}
\put(135,5){\makebox(0,0)[b]{\tiny $k$}}
\put(155,5){\makebox(0,0)[b]{\tiny $k\pl 1$}}
\put(175,5){\makebox(0,0)[b]{\tiny $k\pl 2$}}
\put(195,5){\makebox(0,0)[b]{\tiny $k\pl 3$}}

\put(115,32){\makebox(0,0)[b]{\tiny $1$}}
\put(155,32){\makebox(0,0)[b]{\tiny $k$}}
\put(195,32){\makebox(0,0)[b]{\tiny $k\pl 1$}}

\put(115,10){\line(0,1){20}} \put(155,30){\line(0,-1){10}}
\put(155,20){\line(-2,-1){20}} \put(155,20){\line(2,-1){20}}
\put(195,30){\line(0,-1){20}}

\put(155,13){\makebox(0,0)[b]{\tiny $\beta$}}

\put(155,11){\circle*{2}}

\end{picture}

\begin{picture}(310,40)

\put(22,20){\makebox(0,0)[l]{$c^\alpha_k$}}

\put(125,20){\makebox(0,0)[l]{$\ldots$}}
\put(201,20){\makebox(0,0)[r]{$\ldots$}}

\put(115,5){\makebox(0,0)[b]{\tiny $1$}}
\put(155,5){\makebox(0,0)[b]{\tiny $k$}}
\put(175,5){\makebox(0,0)[b]{\tiny $k\pl 1$}}

\put(115,32){\makebox(0,0)[b]{\tiny $1$}}
\put(155,32){\makebox(0,0)[b]{\tiny $k$}}
\put(175,32){\makebox(0,0)[b]{\tiny $k\pl 1$}}

\put(115,10){\line(0,1){20}} \put(155,30){\line(0,-1){20}}
\put(175,30){\line(0,-1){20}}

\put(154,20){\makebox(0,0)[r]{\tiny $\alpha$}}

\end{picture}
\end{center}
The composition of Frobenius split equivalences is made according
to the following reductions, which are correlated with the
equations of ${\cal L}_\omega$ on the left of the following
pictures, for ${l\leq k}$:
\begin{center}
\begin{picture}(310,60)

\put(-13,30){\makebox(0,0)[l]{$(aa)$}}
\put(32,31){\makebox(0,0)[l]{$a^\alpha_k a^\beta_l=a^\beta_l
a^\alpha_{k+2}$}}

\put(150,54){\makebox(0,0)[b]{\tiny $l$}}
\put(170,54){\makebox(0,0)[b]{\tiny $l\pl 2$}}
\put(190,54){\makebox(0,0)[b]{\tiny $k\pl 2$}}
\put(210,54){\makebox(0,0)[b]{\tiny $k\pl 4$}}

\put(160,5){\makebox(0,0)[b]{\tiny $l$}}
\put(200,5){\makebox(0,0)[b]{\tiny $k$}}

\put(160,10){\line(0,1){20}} \put(160,32){\line(0,1){10}}
\put(160,42){\line(-1,1){10}} \put(160,42){\line(1,1){10}}

\put(200,10){\line(0,1){10}} \put(200,20){\line(-1,1){10}}
\put(200,20){\line(1,1){10}} \put(190,32){\line(0,1){20}}
\put(200,32){\line(0,1){20}} \put(210,32){\line(0,1){20}}

\put(160,50){\circle*{2}} \put(200,29){\circle*{2}}

\put(160.5,48.5){\makebox(0,0)[t]{\tiny $\beta$}}
\put(200,27){\makebox(0,0)[t]{\tiny $\alpha$}}

\put(230,31){\makebox(0,0){$\leadsto$}}

\put(250,54){\makebox(0,0)[b]{\tiny $l$}}
\put(270,54){\makebox(0,0)[b]{\tiny $l\pl 2$}}
\put(290,54){\makebox(0,0)[b]{\tiny $k\pl 2$}}
\put(310,54){\makebox(0,0)[b]{\tiny $k\pl 4$}}

\put(260,5){\makebox(0,0)[b]{\tiny $l$}}
\put(300,5){\makebox(0,0)[b]{\tiny $k$}}

\put(300,10){\line(0,1){20}} \put(300,32){\line(0,1){10}}
\put(300,42){\line(-1,1){10}} \put(300,42){\line(1,1){10}}

\put(260,10){\line(0,1){10}} \put(260,20){\line(-1,1){10}}
\put(260,20){\line(1,1){10}} \put(250,32){\line(0,1){20}}
\put(260,32){\line(0,1){20}} \put(270,32){\line(0,1){20}}

\put(300,50){\circle*{2}} \put(260,29){\circle*{2}}

\put(300,48){\makebox(0,0)[t]{\tiny $\alpha$}}
\put(260.5,27.5){\makebox(0,0)[t]{\tiny $\beta$}}

\put(170,31){\makebox(0,0)[l]{$\ldots$}}
\put(280,31){\makebox(0,0)[l]{$\ldots$}}

\end{picture}

\begin{picture}(310,60)

\put(-13,30){\makebox(0,0)[l]{$(c2)$}}
\put(32,31){\makebox(0,0)[l]{$c^\alpha_k
c^\beta_k=c^{\alpha\sharp\beta}_k$}}

\put(180,54){\makebox(0,0)[b]{\tiny $k$}}
\put(180,5){\makebox(0,0)[b]{\tiny $k$}}
\put(280,54){\makebox(0,0)[b]{\tiny $k$}}
\put(280,5){\makebox(0,0)[b]{\tiny $k$}}

\put(180,10){\line(0,1){20}} \put(180,32){\line(0,1){20}}
\put(280,10){\line(0,1){42}}

\put(182,20){\makebox(0,0)[l]{\tiny $\alpha$}}
\put(182,42){\makebox(0,0)[l]{\tiny $\beta$}}
\put(282,31){\makebox(0,0)[l]{\tiny $\alpha\sharp\beta$}}

\put(230,31){\makebox(0,0){$\leadsto$}}

\end{picture}

\begin{picture}(310,60)

\put(-13,30){\makebox(0,0)[l]{$(cc)$}}
\put(10,31){\makebox(0,0)[l]{{for $l<k,\;\;$} $c^\alpha_k
c^\beta_l=c^\beta_l c^\alpha_k$}}

\put(160,54){\makebox(0,0)[b]{\tiny $l$}}
\put(190,54){\makebox(0,0)[b]{\tiny $k$}}
\put(260,54){\makebox(0,0)[b]{\tiny $l$}}
\put(290,54){\makebox(0,0)[b]{\tiny $k$}}
\put(160,5){\makebox(0,0)[b]{\tiny $l$}}
\put(190,5){\makebox(0,0)[b]{\tiny $k$}}
\put(260,5){\makebox(0,0)[b]{\tiny $l$}}
\put(290,5){\makebox(0,0)[b]{\tiny $k$}}

\put(160,10){\line(0,1){20}} \put(160,32){\line(0,1){20}}
\put(180,10){\line(0,1){20}} \put(180,32){\line(0,1){20}}
\put(190,10){\line(0,1){20}} \put(190,32){\line(0,1){20}}

\put(260,10){\line(0,1){20}} \put(260,32){\line(0,1){20}}
\put(280,10){\line(0,1){20}} \put(280,32){\line(0,1){20}}
\put(290,10){\line(0,1){20}} \put(290,32){\line(0,1){20}}

\put(192,20){\makebox(0,0)[l]{\tiny $\alpha$}}
\put(158,42){\makebox(0,0)[r]{\tiny $\beta$}}

\put(292,42){\makebox(0,0)[l]{\tiny $\alpha$}}
\put(258,20){\makebox(0,0)[r]{\tiny $\beta$}}

\put(230,31){\makebox(0,0){$\leadsto$}}

\put(164,31){\makebox(0,0)[l]{$\ldots$}}
\put(264,31){\makebox(0,0)[l]{$\ldots$}}

\end{picture}

\begin{picture}(310,60)

\put(-13,30){\makebox(0,0)[l]{$(ab\:1)$}}
\put(32,31){\makebox(0,0)[l]{$a^\alpha_l b^\beta_{k+2}=b^\beta_k
a^\alpha_l$}}

\put(150,54){\makebox(0,0)[b]{\tiny $l$}}
\put(170,54){\makebox(0,0)[b]{\tiny $l\pl 2$}}
\put(200,54){\makebox(0,0)[b]{\tiny $k\pl 2$}}

\put(160,5){\makebox(0,0)[b]{\tiny $l$}}
\put(190,5){\makebox(0,0)[b]{\tiny $k$}}
\put(210,5){\makebox(0,0)[b]{\tiny $k\pl 2$}}

\put(160,10){\line(0,1){10}} \put(160,20){\line(-1,1){10}}
\put(160,20){\line(1,1){10}} \put(150,32){\line(0,1){20}}
\put(160,32){\line(0,1){20}} \put(170,32){\line(0,1){20}}

\put(200,52){\line(0,-1){10}} \put(200,42){\line(-1,-1){10}}
\put(200,42){\line(1,-1){10}} \put(190,30){\line(0,-1){20}}
\put(200,30){\line(0,-1){20}} \put(210,30){\line(0,-1){20}}

\put(160,29){\circle*{2}} \put(200,33){\circle*{2}}
\put(160,27){\makebox(0,0)[t]{\tiny $\alpha$}}
\put(200,35){\makebox(0,0)[b]{\tiny $\beta$}}

\put(230,31){\makebox(0,0){$\leadsto$}}

\put(250,54){\makebox(0,0)[b]{\tiny $l$}}
\put(270,54){\makebox(0,0)[b]{\tiny $l\pl 2$}}
\put(300,54){\makebox(0,0)[b]{\tiny $k\pl 2$}}

\put(260,5){\makebox(0,0)[b]{\tiny $l$}}
\put(290,5){\makebox(0,0)[b]{\tiny $k$}}
\put(310,5){\makebox(0,0)[b]{\tiny $k\pl 2$}}

\put(260,10){\line(0,1){20}} \put(260,32){\line(0,1){10}}
\put(260,42){\line(-1,1){10}} \put(260,42){\line(1,1){10}}

\put(300,32){\line(0,1){20}} \put(300,30){\line(0,-1){10}}
\put(300,20){\line(-1,-1){10}} \put(300,20){\line(1,-1){10}}

\put(260,51){\circle*{2}} \put(300,11){\circle*{2}}
\put(260,49){\makebox(0,0)[t]{\tiny $\alpha$}}
\put(300,13){\makebox(0,0)[b]{\tiny $\beta$}}

\put(174,31){\makebox(0,0)[l]{$\ldots$}}
\put(275,31){\makebox(0,0)[l]{$\ldots$}}

\end{picture}

\begin{picture}(310,60)

\put(-13,30){\makebox(0,0)[l]{$(ab\:3.1)$}}
\put(32,31){\makebox(0,0)[l]{$a^\alpha_k b^\beta_{k+1}=c^\beta_k
c^\alpha_{k+1}$}}

\put(160,54){\makebox(0,0)[b]{\tiny $k$}}
\put(180,54){\makebox(0,0)[b]{\tiny $k\pl 1$}}

\put(170,5){\makebox(0,0)[b]{\tiny $k$}}
\put(190,5){\makebox(0,0)[b]{\tiny $k\pl 1$}}

\put(170,10){\line(0,1){10}} \put(170,20){\line(-1,1){10}}
\put(170,20){\line(1,1){10}} \put(190,10){\line(0,1){20}}
\put(160,32){\line(0,1){20}} \put(180,52){\line(0,-1){10}}
\put(180,42){\line(-1,-1){10}} \put(180,42){\line(1,-1){10}}

\put(170,29){\circle*{2}} \put(180,33){\circle*{2}}
\put(170,27){\makebox(0,0)[t]{\tiny $\alpha$}}
\put(180,35){\makebox(0,0)[b]{\tiny $\beta$}}

\put(230,31){\makebox(0,0){$\leadsto$}}

\put(270,54){\makebox(0,0)[b]{\tiny $k$}}
\put(290,54){\makebox(0,0)[b]{\tiny $k\pl 1$}}

\put(270,5){\makebox(0,0)[b]{\tiny $k$}}
\put(290,5){\makebox(0,0)[b]{\tiny $k\pl 1$}}

\put(270,10){\line(0,1){20}} \put(270,32){\line(0,1){20}}

\put(290,10){\line(0,1){20}} \put(290,32){\line(0,1){20}}

\put(269,20){\makebox(0,0)[r]{\tiny $\beta$}}
\put(291,42){\makebox(0,0)[l]{\tiny $\alpha$}}

\end{picture}

\begin{picture}(310,60)

\put(-13,30){\makebox(0,0)[l]{$(ab\:3.2)$}}
\put(32,31){\makebox(0,0)[l]{$a^\alpha_{k+1} b^\beta_k=c^\alpha_k
c^\beta_{k+1}$}}

\put(170,54){\makebox(0,0)[b]{\tiny $k$}}
\put(190,54){\makebox(0,0)[b]{\tiny $k\pl 1$}}

\put(160,5){\makebox(0,0)[b]{\tiny $k$}}
\put(180,5){\makebox(0,0)[b]{\tiny $k\pl 1$}}

\put(180,10){\line(0,1){10}} \put(180,20){\line(-1,1){10}}
\put(180,20){\line(1,1){10}} \put(160,10){\line(0,1){20}}
\put(190,32){\line(0,1){20}} \put(170,52){\line(0,-1){10}}
\put(170,42){\line(-1,-1){10}} \put(170,42){\line(1,-1){10}}

\put(170,33){\circle*{2}} \put(180,29){\circle*{2}}
\put(180,27){\makebox(0,0)[t]{\tiny $\alpha$}}
\put(170,35){\makebox(0,0)[b]{\tiny $\beta$}}

\put(230,31){\makebox(0,0){$\leadsto$}}

\put(270,54){\makebox(0,0)[b]{\tiny $k$}}
\put(290,54){\makebox(0,0)[b]{\tiny $k\pl 1$}}

\put(270,5){\makebox(0,0)[b]{\tiny $k$}}
\put(290,5){\makebox(0,0)[b]{\tiny $k\pl 1$}}

\put(270,10){\line(0,1){20}} \put(270,32){\line(0,1){20}}

\put(290,10){\line(0,1){20}} \put(290,32){\line(0,1){20}}

\put(269,20){\makebox(0,0)[r]{\tiny $\alpha$}}
\put(291,42){\makebox(0,0)[l]{\tiny $\beta$}}

\end{picture}

\begin{picture}(310,60)

\put(-13,30){\makebox(0,0)[l]{$(ab\:3.3)$}}
\put(32,31){\makebox(0,0)[l]{$a^\alpha_k
b^\beta_k=c^{\omega^{\alpha\sharp\beta}}_k$}}

\put(180,54){\makebox(0,0)[b]{\tiny $k$}}

\put(180,5){\makebox(0,0)[b]{\tiny $k$}}

\put(180,10){\line(0,1){10}} \put(180,20){\line(-1,1){10}}
\put(180,20){\line(1,1){10}} \put(180,52){\line(0,-1){10}}
\put(180,42){\line(-1,-1){10}} \put(180,42){\line(1,-1){10}}

\put(180,33){\circle*{2}} \put(180,29){\circle*{2}}
\put(180,27){\makebox(0,0)[t]{\tiny $\alpha$}}
\put(180,35){\makebox(0,0)[b]{\tiny $\beta$}}

\put(230,31){\makebox(0,0){$\leadsto$}}

\put(280,54){\makebox(0,0)[b]{\tiny $k$}}
\put(280,5){\makebox(0,0)[b]{\tiny $k$}}

\put(280,10){\line(0,1){42}}

\put(281,31){\makebox(0,0)[l]{\tiny $\omega^{\alpha\sharp\beta}$}}

\end{picture}

\begin{picture}(310,60)

\put(-13,30){\makebox(0,0)[l]{$(ac\:1)$}}
\put(32,31){\makebox(0,0)[l]{$a^\alpha_k c^\gamma_l=c^\gamma_l
a^\alpha_k$}}

\put(160,54){\makebox(0,0)[b]{\tiny $l$}}
\put(180,54){\makebox(0,0)[b]{\tiny $k$}}
\put(200,54){\makebox(0,0)[b]{\tiny $k\pl 2$}}

\put(160,5){\makebox(0,0)[b]{\tiny $l$}}
\put(190,5){\makebox(0,0)[b]{\tiny $k$}}

\put(160,10){\line(0,1){20}} \put(160,32){\line(0,1){20}}
\put(190,10){\line(0,1){10}} \put(190,20){\line(-1,1){10}}
\put(190,20){\line(1,1){10}} \put(180,32){\line(0,1){20}}
\put(190,32){\line(0,1){20}} \put(200,32){\line(0,1){20}}

\put(190,29){\circle*{2}} \put(190,27){\makebox(0,0)[t]{\tiny
$\alpha$}} \put(159,42){\makebox(0,0)[r]{\tiny $\gamma$}}

\put(164,31){\makebox(0,0)[l]{$\ldots$}}

\put(230,31){\makebox(0,0){$\leadsto$}}

\put(260,54){\makebox(0,0)[b]{\tiny $l$}}
\put(280,54){\makebox(0,0)[b]{\tiny $k$}}
\put(300,54){\makebox(0,0)[b]{\tiny $k\pl 2$}}

\put(260,5){\makebox(0,0)[b]{\tiny $l$}}
\put(290,5){\makebox(0,0)[b]{\tiny $k$}}

\put(260,10){\line(0,1){20}} \put(260,32){\line(0,1){20}}

\put(290,10){\line(0,1){20}} \put(290,32){\line(0,1){10}}
\put(290,42){\line(-1,1){10}} \put(290,42){\line(1,1){10}}

\put(290,51){\circle*{2}} \put(290,49){\makebox(0,0)[t]{\tiny
$\alpha$}} \put(259,20){\makebox(0,0)[r]{\tiny $\gamma$}}

\put(270,31){\makebox(0,0)[l]{$\ldots$}}

\end{picture}

\begin{picture}(310,60)

\put(-13,30){\makebox(0,0)[l]{$(ac\:2)$}}
\put(32,31){\makebox(0,0)[l]{$a^\alpha_l c^\gamma_{k+2}=c^\gamma_k
a^\alpha_l$}}

\put(160,54){\makebox(0,0)[b]{\tiny $l$}}
\put(180,54){\makebox(0,0)[b]{\tiny $l\pl 2$}}
\put(200,54){\makebox(0,0)[b]{\tiny $k\pl 2$}}

\put(170,5){\makebox(0,0)[b]{\tiny $l$}}
\put(200,5){\makebox(0,0)[b]{\tiny $k$}}

\put(200,10){\line(0,1){20}} \put(200,32){\line(0,1){20}}
\put(170,10){\line(0,1){10}} \put(170,20){\line(-1,1){10}}
\put(170,20){\line(1,1){10}} \put(160,32){\line(0,1){20}}
\put(170,32){\line(0,1){20}} \put(180,32){\line(0,1){20}}

\put(170,29){\circle*{2}} \put(170,27){\makebox(0,0)[t]{\tiny
$\alpha$}} \put(201,42){\makebox(0,0)[l]{\tiny $\gamma$}}

\put(184,31){\makebox(0,0)[l]{$\ldots$}}

\put(230,31){\makebox(0,0){$\leadsto$}}

\put(260,54){\makebox(0,0)[b]{\tiny $l$}}
\put(280,54){\makebox(0,0)[b]{\tiny $l\pl 2$}}
\put(300,54){\makebox(0,0)[b]{\tiny $k\pl 2$}}

\put(270,5){\makebox(0,0)[b]{\tiny $l$}}
\put(300,5){\makebox(0,0)[b]{\tiny $k$}}

\put(300,10){\line(0,1){20}} \put(300,32){\line(0,1){20}}

\put(270,10){\line(0,1){20}} \put(270,32){\line(0,1){10}}
\put(270,42){\line(-1,1){10}} \put(270,42){\line(1,1){10}}

\put(270,51){\circle*{2}} \put(270,49){\makebox(0,0)[t]{\tiny
$\alpha$}} \put(301,20){\makebox(0,0)[l]{\tiny $\gamma$}}

\put(280,31){\makebox(0,0)[l]{$\ldots$}}

\end{picture}

\begin{picture}(310,60)

\put(-13,30){\makebox(0,0)[l]{$(ac\:3)$}}
\put(32,31){\makebox(0,0)[l]{$a^\alpha_k c^\gamma_{k+1}=
a^{\alpha\sharp\gamma}_k$}}

\put(170,54){\makebox(0,0)[b]{\tiny $k$}}
\put(190,54){\makebox(0,0)[b]{\tiny $k\pl 2$}}

\put(180,5){\makebox(0,0)[b]{\tiny $k$}}

\put(180,10){\line(0,1){10}} \put(180,20){\line(-1,1){10}}
\put(180,20){\line(1,1){10}} \put(170,32){\line(0,1){20}}
\put(180,32){\line(0,1){20}} \put(190,32){\line(0,1){20}}

\put(180,29){\circle*{2}} \put(180,27){\makebox(0,0)[t]{\tiny
$\alpha$}} \put(179,42){\makebox(0,0)[r]{\tiny $\gamma$}}

\put(230,31){\makebox(0,0){$\leadsto$}}

\put(270,54){\makebox(0,0)[b]{\tiny $k$}}
\put(290,54){\makebox(0,0)[b]{\tiny $k\pl 2$}}

\put(280,5){\makebox(0,0)[b]{\tiny $k$}}

\put(280,10){\line(0,1){22}} \put(280,32){\line(-1,2){10}}
\put(280,32){\line(1,2){10}}

\put(280,51){\circle*{2}} \put(280,49){\makebox(0,0)[t]{\tiny
$\alpha\sharp\gamma$}}

\end{picture}
\end{center}
If we disregard the ordinals, then this is exactly like
composition of split equivalences.

There are moreover reductions corresponding to the equations
$(bb)$, ${(ab\:2)}$, ${(bc\:1)}$, ${(bc\:2)}$ and ${(bc\:3)}$ of
\cite{DP03} (Section on \emph{Normal forms in} ${\cal L}_\omega$),
which are analogous to $(aa)$, ${(ab\:1)}$, ${(ac\:1)}$,
${(ac\:2)}$ and ${(ac\:3)}$. We do not mention here trivial
reductions involving $c^0_k$, which is equal to 1. As a limit
case, where ${l=k}$, of the reduction corresponding to $(aa)$ we
have
\begin{center}
\begin{picture}(140,60)

\put(0,54){\makebox(0,0)[b]{\tiny $l$}}
\put(20.5,54){\makebox(0,0)[b]{\tiny $l\pl 2$}}
\put(41.5,54){\makebox(0,0)[b]{\tiny $l\pl 4$}}

\put(30,5){\makebox(0,0)[b]{\tiny $l$}}

\put(30,10){\line(0,1){10}} \put(30,20){\line(-2,1){20}}
\put(30,20){\line(1,1){10}}

\put(10,32){\line(0,1){10}} \put(10,42){\line(-1,1){10}}
\put(10,42){\line(1,1){10}}

\put(30,32){\line(0,1){20}} \put(40,32){\line(0,1){20}}

\put(30,29){\circle*{2}} \put(10,51){\circle*{2}}
\put(10,49){\makebox(0,0)[t]{\tiny $\beta$}}
\put(30,27){\makebox(0,0)[t]{\tiny $\alpha$}}

\put(70,31){\makebox(0,0){$\leadsto$}}

\put(99,54){\makebox(0,0)[b]{\tiny $l$}}
\put(120.5,54){\makebox(0,0)[b]{\tiny $l\pl 2$}}
\put(140.5,54){\makebox(0,0)[b]{\tiny $l\pl 4$}}

\put(110,5){\makebox(0,0)[b]{\tiny $l$}}

\put(110,10){\line(0,1){10}} \put(110,20){\line(-1,1){10}}
\put(110,20){\line(2,1){20}}

\put(130,32){\line(0,1){10}} \put(130,42){\line(-1,1){10}}
\put(130,42){\line(1,1){10}}

\put(100,32){\line(0,1){20}} \put(110,32){\line(0,1){20}}

\put(110,29){\circle*{2}} \put(130,51){\circle*{2}}
\put(110,27){\makebox(0,0)[t]{\tiny $\beta$}}
\put(130,49){\makebox(0,0)[t]{\tiny $\alpha$}}

\end{picture}
\end{center}
and analogously in other limit cases. The limit case ${l=k}$ of
${(ab\:1)}$ corresponds to one of the Frobenius equations:
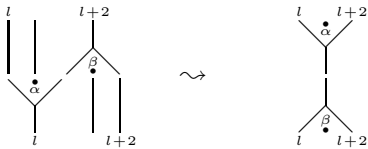
\begin{center}
\begin{picture}(140,60)

\put(0,54){\makebox(0,0)[b]{\tiny $l$}}
\put(32.5,54){\makebox(0,0)[b]{\tiny $l\pl 2$}}

\put(10,5){\makebox(0,0)[b]{\tiny $l$}}
\put(42.5,5){\makebox(0,0)[b]{\tiny $l\pl 2$}}

\put(10,10){\line(0,1){10}} \put(10,20){\line(-1,1){10}}
\put(10,20){\line(1,1){10}}

\put(32,52){\line(0,-1){10}} \put(32,42){\line(-1,-1){10}}
\put(32,42){\line(1,-1){10}}

\put(32,10){\line(0,1){20}} \put(42,10){\line(0,1){20}}
\put(0,32){\line(0,1){20}} \put(10,32){\line(0,1){20}}

\put(10,29){\circle*{2}} \put(32,33){\circle*{2}}
\put(32,35){\makebox(0,0)[b]{\tiny $\beta$}}
\put(10,27){\makebox(0,0)[t]{\tiny $\alpha$}}

\put(70,31){\makebox(0,0){$\leadsto$}}

\put(110,54){\makebox(0,0)[b]{\tiny $l$}}
\put(130,54){\makebox(0,0)[b]{\tiny $l\pl 2$}}

\put(110,5){\makebox(0,0)[b]{\tiny $l$}}
\put(130,5){\makebox(0,0)[b]{\tiny $l\pl 2$}}

\put(120,30){\line(0,-1){10}} \put(120,20){\line(-1,-1){10}}
\put(120,20){\line(1,-1){10}}

\put(120,32){\line(0,1){10}} \put(120,42){\line(-1,1){10}}
\put(120,42){\line(1,1){10}}

\put(120,11){\circle*{2}} \put(120,51){\circle*{2}}
\put(120,13){\makebox(0,0)[b]{\tiny $\beta$}}
\put(120,49){\makebox(0,0)[t]{\tiny $\alpha$}}

\end{picture}
\end{center}

We believe that our Frobenius split equivalences are more handy
than the diagrams that may be found in \cite{KL01} (Appendix~C),
to which they should be equivalent. They are more handy because
the circular forms are coded efficiently by ordinals, while in the
diagrams of \cite{KL01} they make complicated patterns that are
defined in all possible ways in terms of the generators. What
these diagrams miss essentially is the reduction corresponding to
the equation ${(ab\:3.3)}$.

The friezes appropriate for trijunctions (see \cite{DP08b},
Section~8) are such that circular components and circular forms do
not arise. Such friezes can be replaced by maximal split
equivalences, without ordinals. As we said above, in maximal split
equivalences, the odd equivalence classes are completely
determined by the even equivalence classes, and vice versa. By
rejecting the odd equivalence classes, we obtain the split
equivalences that correspond to the categories $S5_{\Box\Diamond}$
and $5S_{\Box\Diamond}$ by the functor $G$; by rejecting the even
equivalence classes, we obtain those that come with the functor
$G^d$ (see \cite{DP08b}, Sections 6-7). Coherence for trijunction
could be proved with respect to nonintersecting split equivalences
for which either odd or even equivalence classes are rejected.

\section{Frobenius monads and matrices}

Let \emph{Mat} be the skeleton of the category $\mbox{\it Vect}_K$
of finite-dimensional vector spaces over the field $K$, with
linear transformations as arrows. The objects of \emph{Mat} are
the natural numbers, which are dimensions of the objects of
$\mbox{\it Vect}_K$, and its arrows are matrices. The category
\emph{Mat} is strictly monoidal (in it the canonical arrows of its
monoidal structure are identity arrows).

In this section we will show how the requirement of having a
faithful functor into \emph{Mat} induces a collapse of the
ordinals of \emph{Frob}. This means that the usual notion of
Frobenius algebra is not exactly caught by the notion of Frobenius
monad. There are further categorial equations implicit in the
notion of Frobenius algebra, which do not hold in every Frobenius
monad. We will describe in this section these equations, and show
their necessity. We leave open the question whether they are also
sufficient to describe categorially the notion of Frobenius
algebra.

There is no faithful monoidal functor from the strictly monoidal
category \emph{Frob} into \emph{Mat}. A necessary condition to
obtain such a functor would be to extend the definition of
\emph{Frob} with some new equations, for whose formulation we need
the following abbreviations:
\begin{tabbing}
\hspace{7.5em}\= $(\delta^\Box_n)^0$\=
$=\mj_{n+1}$,\hspace{9em}\=$(\delta^\Diamond_n)^0$\= $=\mj_{n+1}$,
\\[1.5ex]
\> \hspace{-1em}$(\delta^\Box_n)^{k+1}$\>
$=\delta^\Box_{n+k}\cirk(\delta^\Box_n)^k$,\>
\hspace{-1em}$(\delta^\Diamond_n)^{k+1}$\>
$=(\delta^\Diamond_n)^k\cirk\delta^\Diamond_{n+k}$,
\\[2ex]
\centerline{$\Phi^k_n=_{df}\varepsilon^\Box_n\cirk(\delta^\Diamond_n)^k\cirk(\delta^\Box_n)^k\cirk\varepsilon^\Diamond_n$.}
\end{tabbing}
Our new equations are then all of the following equations, for
${k,n\geq 0}$:
\begin{tabbing}
\hspace{1.7em}$(\Phi)$\hspace{13em}$\Phi^k_n=M^n\Phi^k_0$,
\end{tabbing}
where $M^n$ is a sequence of ${n\geq 0}$ occurrences of $M$.
Equations with the same force as $(\Phi)$, which we will also call
$(\Phi)$, are, for ${k,n\geq 0}$,
\begin{tabbing}
\hspace{1.7em}$(\Phi)$\hspace{13em}\= $\Phi^k_n$ \=
$=M^n\Phi^k_0$,\kill

\> \hspace{-.9em}$\Phi^k_{n+1}$\> $=M\Phi^k_n$.
\end{tabbing}

These equations do not hold in \emph{Frob}, as can be seen with
the help of the monoid ${\cal L}_\omega$, where the corresponding
equations
\begin{tabbing}
\hspace{1.7em}\=$(\Phi)$\hspace{13em}\= $\Phi^k_n$ \=
$=M^n\Phi^k_0$,\kill

\> $(\Phi c)$\> \hspace{-1em}$c^{\omega^k}_{2n+1}$\>
$=c^{\omega^k}_1$
\end{tabbing}
do not hold. These equations hold in the monoid ${\cal K}_\omega$
of \cite{DP03}.

Let the category $\mbox{\it Frob}'$ be defined like \emph{Frob}
save that we have in addition all the equations $(\Phi)$, and let
${\cal L}_\omega'$ be the monoid defined like ${\cal L}_\omega$
save that we have in addition all the equations $(\Phi c)$. If all
the subscripts $n$ that may be found in defining $\Phi^k_n$ are
replaced by $A$, while ${n\pl 1}$ and ${n\pl k}$ are replaced
respectively by $MA$ and ${M^kA}$, then the equations $(\Phi)$
become
\begin{tabbing}
\hspace{1.7em}\=$(\Phi)$\hspace{8em}\= $\Phi^k_n$ \=
$=M^n\Phi^k_0$,\kill

\>\> \hspace{-1.2em}$\Phi^k_{M^{\!n}\!A}$\> $=M^n\Phi^k_A$\quad
or\quad $\Phi^k_{MA}=M\Phi^k_A$,
\end{tabbing}
which we will also call $(\Phi)$, and which are the equations
characterizing the class of Frobenius monads in which $\mbox{\it
Frob}'$ is the free one generated with a single object.

In the language of the free self-adjunction of Section~4, let
$\kappa^0_{2n+1}$ stand for $\mj_{2n+1}$, and let
$\kappa^{k+1}_{2n+1}$ be
${\kappa^k_{2n+1}\cirk\varphi_{2n+1}\cirk\gamma_{2n+1}}$. Consider
then the category $\cal S'$ constructed like the category $\cal S$
of the free self-adjunction save that we have in addition for
every ${k,n\geq 0}$ the equation
\[
\varphi_{2n}\cirk
F\kappa^k_{2n+1}\cirk\gamma_{2n}=F^{2n}(\varphi_0\cirk
F\kappa^k_1\cirk\gamma_0),
\]
where $F^m$ is a sequence of ${m\geq 0}$ occurrences of $F$. The
category $\cal S'$ is related to $\mbox{\it Frob}'$ as the
category $\cal S$ is related to \emph{Frob}; this is shown as in
Section~5. On the other hand, $\cal S'$ is related to ${\cal
L}_\omega'$ as $\cal S$ is related to ${\cal L}_\omega$; this is
shown as in \cite{DP03} (Section on ${\cal L}_c$ \emph{and} ${\cal
L}_\omega$).

We can infer that $\mbox{\it Frob}'$ is isomorphic to a category
whose arrows are the elements of the monoid ${\cal L}_\omega'$
with types associated to them (see the preceding section). This
result may be understood as a coherence result, which provides a
decision procedure for equality of arrows in $\mbox{\it Frob}'$.
The normal form involved in this decision procedure would serve
also for the isomorphism with the category $\mbox{\it Frz}'$,
which we will consider in a moment. We will deal with this normal
form later (see the second paragraph after the proof of Lemma
$2m\pl 2$).

One could consider a category $\mbox{\it Frz}'$ analogous to the
category \emph{Frz} of the previous section, which would be
isomorphic to our category derived from ${\cal L}_\omega'$. We
will not describe $\mbox{\it Frz}'$ in detail, but just make a few
indications. For the arrows of $\mbox{\it Frz}'$ we would take,
instead of friezes, two-manifolds made out of friezes in the
following way. The regions of friezes may be chessboard-coloured
by making the leftmost region white, and then alternating black
and white for subsequent regions. For example, one of the friezes
we had above is chessboard-coloured as follows:
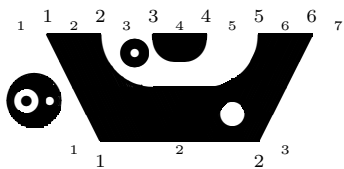
\begin{center}
\begin{picture}(100,60)

\put(0,55){\makebox(0,0)[b]{\scriptsize $1$}}
\put(20,55){\makebox(0,0)[b]{\scriptsize $2$}}
\put(40,55){\makebox(0,0)[b]{\scriptsize $3$}}
\put(60,55){\makebox(0,0)[b]{\scriptsize $4$}}
\put(80,55){\makebox(0,0)[b]{\scriptsize $5$}}
\put(100,55){\makebox(0,0)[b]{\scriptsize $6$}}

\put(20,0){\makebox(0,0)[b]{\scriptsize $1$}}
\put(80,0){\makebox(0,0)[b]{\scriptsize $2$}}

\put(-10,52){\makebox(0,0)[b]{\tiny $1$}}
\put(10,52){\makebox(0,0)[b]{\tiny $2$}}
\put(30,52){\makebox(0,0)[b]{\tiny $3$}}
\put(50,52){\makebox(0,0)[b]{\tiny $4$}}
\put(70,52){\makebox(0,0)[b]{\tiny $5$}}
\put(90,52){\makebox(0,0)[b]{\tiny $6$}}
\put(110,52){\makebox(0,0)[b]{\tiny $7$}}

\put(10,8){\makebox(0,0)[t]{\tiny $1$}}
\put(50,8){\makebox(0,0)[t]{\tiny $2$}}
\put(90,8){\makebox(0,0)[t]{\tiny $3$}}

\put(40,50){\line(1,0){20}}

\thicklines

\put(40,50){\line(1,0){20}}

\put(0,50){\line(1,0){20}} \put(80,50){\line(1,0){20}}
\put(20,10){\line(-1,2){20}} \put(17,49.2){\line(1,0){3}}
\put(17,48.4){\line(1,0){3}} \put(17,47.6){\line(1,0){3}}
\put(17,46.8){\line(1,0){3}} \put(17,46){\line(1,0){3}}
\put(80,49.2){\line(1,0){2}} \put(80,48.4){\line(1,0){2}}
\put(80,47.6){\line(1,0){2}} \put(80,46.8){\line(1,0){2}}
\put(80,46){\line(1,0){2}} \put(79.8,45.2){\line(1,0){2}}
\put(74,25){\line(1,0){10}}

\multiput(59.7,30)(0,-.6){9}{\line(1,0){13}}
\multiput(67.4,24.4)(-.3,-.4){14}{\line(-1,0){7}}
\multiput(74.5,30)(.75,.4){16}{\line(-1,0){10}}
\multiput(80,24.4)(.2,-.4){18}{\line(-1,0){7}}
\multiput(63,10.4)(0,.4){12}{\line(1,0){11}}
\put(20,10){\line(1,0){60}}

\multiput(20,10)(.4,0){44}{\line(-1,2){20}}
\put(37.8,10){\line(-1,2){14.5}} \put(38.2,10){\line(-1,2){14}}
\put(38.6,10){\line(-1,2){13}} \put(39,10){\line(-1,2){12.5}}
\put(39.4,10){\line(-1,2){12}} \put(39.8,10){\line(-1,2){11.9}}
\put(40.2,10){\line(-1,2){11.8}}
\multiput(41,10)(.8,0){2}{\line(-1,2){11.15}}
\put(42.6,10){\line(-1,2){10.5}}
\multiput(43.4,10)(.8,0){4}{\line(-1,2){10.3}}
\multiput(37,10)(.4,0){80}{\line(-1,2){9.8}}

\put(80,10){\line(1,2){20}}
\multiput(80,10)(-.4,0){24}{\line(1,2){20}}
\multiput(58.5,10)(-.4,0){40}{\line(1,2){9.8}}
\multiput(79,25)(-.4,0){24}{\line(1,2){12.5}}

\put(50,50){\oval(20,20)[b]} \put(50,50){\oval(19.5,19.5)[b]}
\put(50,50){\oval(19,19)[b]} \put(50,50){\oval(18.5,18.5)[b]}
\put(50,50){\oval(18,18)[b]} \put(50,50){\oval(17.5,17.5)[b]}
\put(50,50){\oval(17,17)[b]} \put(50,50){\oval(16.5,16.5)[b]}
\put(50,50){\oval(16,16)[b]} \put(50,50){\oval(15.5,15.5)[b]}
\put(50,50){\oval(15,15)[b]} \put(50,50){\oval(14.5,14.5)[b]}
\put(50,50){\oval(14,14)[b]} \put(50,50){\oval(13.5,13.5)[b]}
\put(50,50){\oval(13,13)[b]} \put(50,50){\oval(12.5,12.5)[b]}
\put(50,50){\oval(12,12)[b]} \put(50,50){\oval(11.5,11.5)[b]}
\put(50,50){\oval(11,11)[b]} \put(50,50){\oval(10.5,10.5)[b]}
\put(50,50){\oval(10,10)[b]} \put(50,50){\oval(9.5,9.5)[b]}
\put(50,50){\oval(9,9)[b]} \put(50,50){\oval(8.5,8.5)[b]}
\put(50,50){\oval(8,8)[b]} \put(50,50){\oval(7.5,7.5)[b]}
\put(50,50){\oval(7,7)[b]} \put(50,50){\oval(6.5,6.5)[b]}
\put(50,50){\oval(6,6)[b]} \put(50,50){\oval(5.5,5.5)[b]}
\put(50,50){\oval(5,5)[b]} \put(50,50){\oval(4.5,4.5)[b]}
\put(50,50){\oval(4,4)[b]} \put(50,50){\oval(3.5,3.5)[b]}
\put(50,50){\oval(3,3)[b]} \put(50,50){\oval(2.5,2.5)[b]}
\put(50,50){\oval(2,2)[b]} \put(50,50){\oval(1.5,1.5)[b]}
\put(50,50){\oval(1,1)[b]}

\put(50,50){\oval(60,40)[b]} \put(50,50){\oval(61,41)[b]}

\put(33,43){\circle{10}} \put(33,43){\circle{9.5}}
\put(33,43){\circle{9}} \put(33,43){\circle{8.5}}
\put(33,43){\circle{8}} \put(33,43){\circle{7.5}}
\put(33,43){\circle{7}} \put(33,43){\circle{6.5}}
\put(33,43){\circle{6}} \put(33,43){\circle{5.5}}
\put(33,43){\circle{5}} \put(33,43){\circle{4.5}}
\put(33,43){\circle{4}}

\put(70,20){\circle{10}} \put(70,20){\circle{11}}
\put(70,20){\circle{12}} \put(70,20){\circle{13}}
\put(-5,25){\circle{20}} \put(1,25){\circle{4}}
\put(1,25){\circle{5}} \put(1,25){\circle{6}}
\put(1,25){\circle{7}} \put(-8,25){\circle{10}}
\put(-8,25){\circle{11}} \put(-8,25){\circle{12}}
\put(-8,25){\circle{13}} \put(-8,25){\circle{14}}
\put(-8,25){\circle*{4}}

\put(-7,34.6){\line(1,0){4}} \put(-7,34.2){\line(1,0){4}}
\put(-9.8,33.8){\line(1,0){9.6}} \put(-9.8,33.4){\line(1,0){9}}
\put(-11,33){\line(1,0){12}} \put(-11,32.6){\line(1,0){11}}
\put(-12,32.2){\line(1,0){14}} \put(-12,31.8){\line(1,0){13}}
\put(-13,31.4){\line(1,0){15}} \put(-13,31){\line(1,0){15}}
\put(-13.5,30.6){\line(1,0){16}} \put(-13.5,30.2){\line(1,0){16}}
\put(-14,29.8){\line(1,0){5}} \put(-14,29.4){\line(1,0){3.5}}
\put(-14.3,29){\line(1,0){3.5}} \put(-14.3,21){\line(1,0){3.5}}
\put(-14,20.6){\line(1,0){3.5}} \put(-14,20.2){\line(1,0){5}}
\put(-13.5,19.8){\line(1,0){16}} \put(-13.5,19.4){\line(1,0){17}}
\put(-13,19){\line(1,0){15}} \put(-13,18.6){\line(1,0){16}}
\put(-12,18.2){\line(1,0){13}} \put(-12,17.8){\line(1,0){14}}
\put(-11,17.4){\line(1,0){11}} \put(-11,17){\line(1,0){12}}
\put(-9.8,16.6){\line(1,0){9}} \put(-9.8,16.2){\line(1,0){9.6}}
\put(-8,15.8){\line(1,0){6}} \put(-8,15.4){\line(1,0){6.5}}
\put(-4,30.2){\line(1,0){7}} \put(-4,29.8){\line(1,0){7.5}}
\put(-3.5,29.4){\line(1,0){7}} \put(-3.5,29){\line(1,0){7.5}}
\put(-3.2,28.6){\line(1,0){7}} \put(-3.2,28.2){\line(1,0){7}}
\put(-3.2,27.8){\line(1,0){8}} \put(-3.2,27.4){\line(1,0){8}}
\put(-4,20.6){\line(1,0){8}} \put(-4,20.2){\line(1,0){8}}
\put(-3.5,21.4){\line(1,0){8}} \put(-3.5,21){\line(1,0){8}}
\put(-3.2,22.2){\line(1,0){8}} \put(-3.2,21.8){\line(1,0){8}}
\put(-3.2,22.6){\line(1,0){8}}

\end{picture}
\end{center}
Then consider the two-manifolds with boundary made of the compact
black regions, which we will call \emph{black friezes}, and on
black friezes consider the equivalence relation based on
homeomorphisms that preserve all the points on the top and bottom
line (this is like the $\cal K$-equivalence of \cite{DP03},
Section on \emph{Friezes}). So the following black frieze would be
equivalent to the black frieze above:
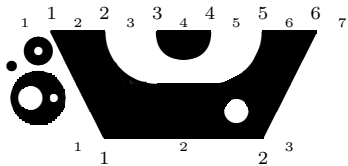
\begin{center}
\begin{picture}(100,60)

\put(0,55){\makebox(0,0)[b]{\scriptsize $1$}}
\put(20,55){\makebox(0,0)[b]{\scriptsize $2$}}
\put(40,55){\makebox(0,0)[b]{\scriptsize $3$}}
\put(60,55){\makebox(0,0)[b]{\scriptsize $4$}}
\put(80,55){\makebox(0,0)[b]{\scriptsize $5$}}
\put(100,55){\makebox(0,0)[b]{\scriptsize $6$}}

\put(20,0){\makebox(0,0)[b]{\scriptsize $1$}}
\put(80,0){\makebox(0,0)[b]{\scriptsize $2$}}

\put(-10,52){\makebox(0,0)[b]{\tiny $1$}}
\put(10,52){\makebox(0,0)[b]{\tiny $2$}}
\put(30,52){\makebox(0,0)[b]{\tiny $3$}}
\put(50,52){\makebox(0,0)[b]{\tiny $4$}}
\put(70,52){\makebox(0,0)[b]{\tiny $5$}}
\put(90,52){\makebox(0,0)[b]{\tiny $6$}}
\put(110,52){\makebox(0,0)[b]{\tiny $7$}}

\put(10,8){\makebox(0,0)[t]{\tiny $1$}}
\put(50,8){\makebox(0,0)[t]{\tiny $2$}}
\put(90,8){\makebox(0,0)[t]{\tiny $3$}}

\put(40,50){\line(1,0){20}}

\thicklines

\put(40,50){\line(1,0){20}}

\put(0,50){\line(1,0){20}} \put(80,50){\line(1,0){20}}
\put(20,10){\line(-1,2){20}} \put(17,49.2){\line(1,0){3}}
\put(17,48.4){\line(1,0){3}} \put(17,47.6){\line(1,0){3}}
\put(17,46.8){\line(1,0){3}} \put(17,46){\line(1,0){3}}
\put(80,49.2){\line(1,0){2}} \put(80,48.4){\line(1,0){2}}
\put(80,47.6){\line(1,0){2}} \put(80,46.8){\line(1,0){2}}
\put(80,46){\line(1,0){2}} \put(79.8,45.2){\line(1,0){2}}
\put(74,25){\line(1,0){10}}

\multiput(59.7,30)(0,-.6){9}{\line(1,0){13}}
\multiput(67.4,24.4)(-.3,-.4){14}{\line(-1,0){7}}
\multiput(74.5,30)(.75,.4){16}{\line(-1,0){10}}
\multiput(80,24.4)(.2,-.4){18}{\line(-1,0){7}}
\multiput(63,10.4)(0,.4){12}{\line(1,0){11}}
\put(20,10){\line(1,0){60}}

\multiput(20,10)(.4,0){44}{\line(-1,2){20}}
\put(37.8,10){\line(-1,2){14.5}} \put(38.2,10){\line(-1,2){14}}
\put(38.6,10){\line(-1,2){13}} \put(39,10){\line(-1,2){12.5}}
\put(39.4,10){\line(-1,2){12}} \put(39.8,10){\line(-1,2){11.9}}
\put(40.2,10){\line(-1,2){11.8}}
\multiput(41,10)(.8,0){2}{\line(-1,2){11.15}}
\put(42.6,10){\line(-1,2){10.5}}
\multiput(43.4,10)(.8,0){4}{\line(-1,2){10.3}}
\multiput(37,10)(.4,0){80}{\line(-1,2){9.8}}

\put(80,10){\line(1,2){20}}
\multiput(80,10)(-.4,0){24}{\line(1,2){20}}
\multiput(58.5,10)(-.4,0){40}{\line(1,2){9.8}}
\multiput(79,25)(-.4,0){24}{\line(1,2){12.5}}

\put(50,50){\oval(20,20)[b]} \put(50,50){\oval(19.5,19.5)[b]}
\put(50,50){\oval(19,19)[b]} \put(50,50){\oval(18.5,18.5)[b]}
\put(50,50){\oval(18,18)[b]} \put(50,50){\oval(17.5,17.5)[b]}
\put(50,50){\oval(17,17)[b]} \put(50,50){\oval(16.5,16.5)[b]}
\put(50,50){\oval(16,16)[b]} \put(50,50){\oval(15.5,15.5)[b]}
\put(50,50){\oval(15,15)[b]} \put(50,50){\oval(14.5,14.5)[b]}
\put(50,50){\oval(14,14)[b]} \put(50,50){\oval(13.5,13.5)[b]}
\put(50,50){\oval(13,13)[b]} \put(50,50){\oval(12.5,12.5)[b]}
\put(50,50){\oval(12,12)[b]} \put(50,50){\oval(11.5,11.5)[b]}
\put(50,50){\oval(11,11)[b]} \put(50,50){\oval(10.5,10.5)[b]}
\put(50,50){\oval(10,10)[b]} \put(50,50){\oval(9.5,9.5)[b]}
\put(50,50){\oval(9,9)[b]} \put(50,50){\oval(8.5,8.5)[b]}
\put(50,50){\oval(8,8)[b]} \put(50,50){\oval(7.5,7.5)[b]}
\put(50,50){\oval(7,7)[b]} \put(50,50){\oval(6.5,6.5)[b]}
\put(50,50){\oval(6,6)[b]} \put(50,50){\oval(5.5,5.5)[b]}
\put(50,50){\oval(5,5)[b]} \put(50,50){\oval(4.5,4.5)[b]}
\put(50,50){\oval(4,4)[b]} \put(50,50){\oval(3.5,3.5)[b]}
\put(50,50){\oval(3,3)[b]} \put(50,50){\oval(2.5,2.5)[b]}
\put(50,50){\oval(2,2)[b]} \put(50,50){\oval(1.5,1.5)[b]}
\put(50,50){\oval(1,1)[b]}

\put(50,50){\oval(60,40)[b]} \put(50,50){\oval(61,41)[b]}

\put(-5,43){\circle{10}} \put(-5,43){\circle{9.5}}
\put(-5,43){\circle{9}} \put(-5,43){\circle{8.5}}
\put(-5,43){\circle{8}} \put(-5,43){\circle{7.5}}
\put(-5,43){\circle{7}} \put(-5,43){\circle{6.5}}
\put(-5,43){\circle{6}} \put(-5,43){\circle{5.5}}
\put(-5,43){\circle{5}} \put(-5,43){\circle{4.5}}
\put(-5,43){\circle{4}}

\put(70,20){\circle{10}} \put(70,20){\circle{11}}
\put(70,20){\circle{12}} \put(70,20){\circle{13}}
\put(-5,25){\circle{20}} \put(1,25){\circle{4}}
\put(1,25){\circle{5}} \put(1,25){\circle{6}}
\put(1,25){\circle{7}} \put(-8,25){\circle{10}}
\put(-8,25){\circle{11}} \put(-8,25){\circle{12}}
\put(-8,25){\circle{13}} \put(-8,25){\circle{14}}
\put(-15,37){\circle*{4}}

\put(-7,34.6){\line(1,0){4}} \put(-7,34.2){\line(1,0){4}}
\put(-9.8,33.8){\line(1,0){9.6}} \put(-9.8,33.4){\line(1,0){9}}
\put(-11,33){\line(1,0){12}} \put(-11,32.6){\line(1,0){11}}
\put(-12,32.2){\line(1,0){14}} \put(-12,31.8){\line(1,0){13}}
\put(-13,31.4){\line(1,0){15}} \put(-13,31){\line(1,0){15}}
\put(-13.5,30.6){\line(1,0){16}} \put(-13.5,30.2){\line(1,0){16}}
\put(-14,29.8){\line(1,0){5}} \put(-14,29.4){\line(1,0){3.5}}
\put(-14.3,29){\line(1,0){3.5}} \put(-14.3,21){\line(1,0){3.5}}
\put(-14,20.6){\line(1,0){3.5}} \put(-14,20.2){\line(1,0){5}}
\put(-13.5,19.8){\line(1,0){16}} \put(-13.5,19.4){\line(1,0){17}}
\put(-13,19){\line(1,0){15}} \put(-13,18.6){\line(1,0){16}}
\put(-12,18.2){\line(1,0){13}} \put(-12,17.8){\line(1,0){14}}
\put(-11,17.4){\line(1,0){11}} \put(-11,17){\line(1,0){12}}
\put(-9.8,16.6){\line(1,0){9}} \put(-9.8,16.2){\line(1,0){9.6}}
\put(-8,15.8){\line(1,0){6}} \put(-8,15.4){\line(1,0){6.5}}
\put(-4,30.2){\line(1,0){7}} \put(-4,29.8){\line(1,0){7.5}}
\put(-3.5,29.4){\line(1,0){7}} \put(-3.5,29){\line(1,0){7.5}}
\put(-3.2,28.6){\line(1,0){7}} \put(-3.2,28.2){\line(1,0){7}}
\put(-3.2,27.8){\line(1,0){8}} \put(-3.2,27.4){\line(1,0){8}}
\put(-4,20.6){\line(1,0){8}} \put(-4,20.2){\line(1,0){8}}
\put(-3.5,21.4){\line(1,0){8}} \put(-3.5,21){\line(1,0){8}}
\put(-3.2,22.2){\line(1,0){8}} \put(-3.2,21.8){\line(1,0){8}}
\put(-3.2,22.6){\line(1,0){8}}

\end{picture}
\end{center}
The category $\mbox{\it Frz}'$ is related to the category
\emph{2Cob} of \cite{K03} (Section 1.4), whose arrows are
cobordisms of dimension 2. An arrow of $\mbox{\it Frz}'$ may be
conceived as a kind of ``thin'' cobordism.

As we associated the category \emph{Frobse} to \emph{Frz}, so we
may look for a category $\mbox{\it Frobse}'$ like \emph{Frobse} to
associate to $\mbox{\it Frz}'$. We will previously demonstrate
however the necessity of the equations $(\Phi)$ for faithful
monoidal functors into \emph{Mat}, and consider the consequences
for ordinals of having $(\Phi)$ and related equations.

The necessity of $(\Phi)$ follows from the fact that \emph{Mat} is
a symmetric strictly monoidal category, which has a symmetry
natural isomorphism $c_{n,m}\!:n\otimes m\str m\otimes n$ for
which we have the equation
\begin{tabbing}
\hspace{1.7em}${(c1)}$\hspace{11em}$c_{1,m}=c_{m,1}=\mj_m$
\end{tabbing}
(where 1 in the subscripts of $c$ is the unit object of
\emph{Mat}). Hence, for every arrow ${f\!:1\str 1}$ of \emph{Mat},
we have
\[
\mj_m\otimes f=(\mj_m\otimes f)\cirk c_{1,m}=
c_{1,m}\cirk(f\otimes\mj_m)=f\otimes\mj_m.
\]
Since for every monoidal functor $G$ from \emph{Frob} to
\emph{Mat} we have $G0=1$ (where 0 is the unit object of
\emph{Frob}), and since $G\Phi^k_n$ is of the form ${\mj_{n\cdot
p}\otimes f}$ for ${f\!:1\str 1}$, we have $G\Phi^k_n=G\Phi^k_0$.
So, from the faithfulness of $G$, the equation $(\Phi)$ follows.

In the reasoning above $c$ can be a braiding natural isomorphism,
instead of a symmetry natural isomorphism. We would have the
equation ${(c1)}$, and the equation $(\Phi)$ would again be
imposed by the faithfulness of $G$. So we could replace \emph{Mat}
by a braided strictly monoidal category (cf.\ \cite{K03}, Section
3.6.27).

We defined above the monoid ${\cal L}_\omega'$ as ${\cal
L}_\omega$ with the equation $(\Phi c)$ added. In ${\cal
L}_\omega'$ the hierarchy of $\varepsilon_0$ collapses to
$\omega^\omega$. This means that every element of ${\cal
L}_\omega'$ is definable in terms of $e^\beta_n$, for $e$ being
$a$, $b$ or $c$, and ${\beta\in\omega^\omega}$. We can restrict
the terms $e^\beta_n$ even further, to those in the following
table, without altering the structure of the normal form for
${\cal L}_\omega$ of \cite{DP03} (Section on \emph{Normal forms
in} ${\cal L}_\omega$):
\begin{center}
\begin{tabular}{c|c|c}
$e$ & $n$ & $\beta$
\\[.5ex]
\hline $c$ & 1 & $\hspace{0.7em}\beta\in\omega^\omega$
\\[.3ex]
$c$ & $2m\pl 2$ & $\beta\in\omega$
\\[.3ex]
$a$ and $b$ & $2m\pl 1$ & $\beta\in\omega$
\\[.3ex]
$a$ and $b$ & $2m\pl 2$ & $\beta=0$
\end{tabular}
\end{center}
This is shown as follows.

By Cantor's Normal Form Theorem (see, for example, \cite{KM},
VII.7, Theorem~2, p.\ 248, or \cite{L79}, IV.2, Theorem 2.14, p.\
127), for every ordinal ${\alpha>0}$ in $\varepsilon_0$ there is a
unique finite ordinal ${n\geq 1}$ and a unique sequence of
ordinals $\alpha_1\geq\ldots\geq\alpha_n$ contained in $\alpha$,
i.e.\ lesser than $\alpha$, such that
$\alpha=\omega^{\alpha_1}\sharp\ldots\sharp\,\omega^{\alpha_n}$.
So every ordinal in $\varepsilon_0$ can be named by using the
operations of the monoid
${\langle\varepsilon_0,\sharp,0,\omega^{-}\rangle}$ mentioned in
the previous section.

Let $\beta_0$ be $\omega^0$, which is equal to 1, and let
${\beta_k\!:\varepsilon_0^k\str\varepsilon_0}$, for ${k\geq 1}$,
be defined by
\begin{tabbing}
\hspace{9.7em}\= If $k=0$, then $\beta_0'$ \=
$=\omega^0=1=\beta_0$.\kill

\> \hspace{1em}$\beta_k(\alpha_1,\ldots,\alpha_k)$ \>
$=\omega^{\omega^{\alpha_1}\sharp\ldots\sharp\,\omega^{\alpha_k}}$.
\end{tabbing}
By Cantor's Normal Form Theorem, to name the ordinals in
$\varepsilon_0$ we can replace the unary operation $\omega^{-}$ by
the operations $\beta_k$ for every ${k\geq 0}$. So the name of
every ordinal in $\varepsilon_0$ can be written in terms of $0$,
$\sharp$ and $\beta_k$. We proceed by induction on the complexity
of such a name to define the map $'$ from $\varepsilon_0$ to
$\omega^\omega$:
\begin{tabbing}
\hspace{9.2em}\= If $k=0$, then $\beta_0'$ \=
$=\omega^0=1=\beta_0$.\kill

\>\hspace{6.3em} $0'$ \> $=0$,
\\[1.5ex]
\> \hspace{3.6em}$(\alpha_1\sharp\,\alpha_2)'$ \>
$=\alpha_1'\sharp\,\alpha_2'$,
\\[1.5ex]
\>\hspace{6.4em}$\beta_0'$ \> $=\omega^0=1=\beta_0$,
\\*[1.5ex]
\>\hspace{.9em}$\beta_k(\alpha_1,\ldots,\alpha_k)'$ \>
$=\omega^k\sharp\,\alpha_1'\sharp\ldots\sharp\,\alpha_k'$,\hspace{1em}
for $k\geq 1$.
\end{tabbing}
We can then prove the following lemmata.

\prop{Lemma $2m\pl 1$}{In ${\cal L}_\omega'$, for every $m\geq 0$,
we have $c^\alpha_{2m+1}=c^{\alpha'}_1$.}

\dkz We proceed by induction on the size of $\alpha$. If
${\alpha=0}$, then we use the following equation of ${\cal
L}_\omega$:
\begin{tabbing}
\hspace{3.5em}$c^0_{2m+1}=c^0_1=\mj$.
\end{tabbing}
In the induction step we have
\begin{tabbing}
\hspace{3.5em}\=$c^{\alpha_1\sharp\,\alpha_2}_{2m+1}$ \=
$=c^{\alpha_1'\sharp\,\alpha_2'}_{2m+1}$, \hspace{1em}by ${(c2)}$
and the induction hypothesis,
\\[1.5ex]
\> $c^{\beta_0}_{2m+1}$ \> $=c^{\beta_0'}_1$, \hspace{1em}by
$(\Phi c)$,
\\[2ex]
for $k\geq 1$,\\*[1ex]
\>$c^{\beta_k(\alpha_1,\ldots,\alpha_k)}_{2m+1}=a^0_{2m+1}c^{\omega^{\alpha_1}}_{2m+2}\ldots
c^{\omega^{\alpha_k}}_{2m+2} b^0_{2m+1}$, by ${(ab\:3.3)}$,
${(ac\:3)}$ and ${(c2)}$.
\end{tabbing}
For every $i\in\{1\ldots,k\}$, we have, by the same equations,
\begin{tabbing}
\hspace{3.5em}\=$c^{\alpha_1\sharp\,\alpha_2}_{2m+1}$ \= \kill \>
$c^{\omega^{\alpha_i}}_{2m+2}$ \>
$=a^0_{2m+2}c^{\alpha_i}_{2m+3}b^0_{2m+2}$.
\end{tabbing}
Then, by the induction hypothesis and the equations ${(ac\:1)}$,
${(bc\:1)}$ and ${(c2)}$, for $d^0$ being $\mj$, and $d^{n+1}$
being $d^na^0_{2m+2}b^0_{2m+2}$, we obtain
\begin{tabbing}
\hspace{3.5em}$c^{\beta_k(\alpha_1,\ldots,\alpha_k)}_{2m+1}$ \=
$=a^0_{2m+1}d^k
b^0_{2m+1}c_1^{\alpha_1'\sharp\ldots\sharp\,\alpha_k'}$
\\[1.5ex]
\>
$=c^{\omega^k}_{2m+1}c^{\alpha_1'\sharp\ldots\sharp\,\alpha_k'}_1$,
\hspace{1em}by ${(ab\:3.3)}$ and ${(ac\:3)}$,
\\[1.5ex]
\> $=c^{\beta_k(\alpha_1,\ldots,\alpha_k)'}_1$, \hspace{1em}by
$(\Phi c)$ and ${(c2)}$. \` $\dashv$
\end{tabbing}

\prop{Lemma $2m\pl 2$}{In ${\cal L}_\omega'$, for every $m\geq 0$,
we have $c^{\omega^\alpha}_{2m+2}=c^{\alpha'}_1c^1_{2m+2}$.}

\dkz We have
\begin{tabbing}
\hspace{3.5em}\=$c^{\alpha_1\sharp\,\alpha_2}_{2m+1}$ \= \kill \>
$c^{\omega^{\alpha}}_{2m+2}$ \>
$=a^0_{2m+2}c^{\alpha}_{2m+3}b^0_{2m+2}$, \hspace{1em}by
${(ab\:3.3)}$ and ${(ac\:3)}$,
\\[1.5ex]
\>\> $=c^{\alpha'}_1c^1_{2m+2}$, \hspace{1em}by the preceding
lemma, ${(ac\:1)}$ and ${(ab\:3.3)}$.\` $\dashv$
\end{tabbing}

\vspace{2ex}

With these two lemmata, we can show that the terms $e^\beta_n$ in
the table above are sufficient to define every element of ${\cal
L}_\omega'$ without altering the structure of our normal form.
This is clear for the terms $c^\alpha_n$. We also have
\begin{tabbing}
\hspace{5.6em}\=$c^{\alpha_1\sharp\,\alpha_2}_{2m+1}$ \= \kill

\> $a^\alpha_{2m+2}$ \> $=a^0_{2m+2}c^\alpha_{2m+3}$,
\hspace{1em}by ${(ac\:3)}$,
\\[1.5ex]
\>\> $=c^{\alpha'}_1 a^0_{2m+2}$, \hspace{1em}by Lemma $2m+1$ and
${(ac\:1)}$;
\\[2ex]
\>
\hspace{-2em}$a^{\omega^{\alpha_1}\sharp\ldots\sharp\,\omega^{\alpha_n}}_{2m+1}$
\> $=a^0_{2m+1}c^{\omega^{\alpha_1}}_{2m+2}\ldots
c^{\omega^{\alpha_n}}_{2m+2}$, \hspace{1em}by ${(ac\:3)}$ and
$(c2)$,
\\[1.5ex]
\>\> $=c^{\alpha_1'\sharp\ldots\sharp\,\alpha_n'}a^n_{2m+1}$,
\hspace{1em}by Lemma $2m+2$, ${(ac\:1)}$ and ${(c2)}$,
\end{tabbing}
and analogous equations with $a$ replaced by~$b$.

Consider terms of ${\cal L}_\omega'$ in the form exactly like the
normal form of ${\cal L}_\omega$ in \cite{DP03} save that all the
generators $a^\alpha_i$, $b^\beta_j$ and $c^\gamma_k$ are terms
from our table. We say that such terms are in \emph{normal form}.
This is the normal form we mentioned previously, which we can use
to decide equations in ${\cal L}_\omega'$, and to prove the
isomorphism with $\mbox{\it Frz}'$, along the lines
of~\cite{DP03}.

We can now sketch how the category $\mbox{\it Frobse}'$ analogous
to \emph{Frobse} and isomorphic to $\mbox{\it Frz}'$ would look
like. Its arrows will be based on Frobenius split equivalences
where the function assigning ordinals will follow restrictions in
accordance with our table:

\vspace{1ex}

(1) an even class is mapped to an ordinal in $\omega$,

\vspace{.5ex}

(2) an odd class containing 1 is mapped to an ordinal in
$\omega^\omega$,

\vspace{.5ex}

(3) an odd class not containing 1 is mapped to 0.

\vspace{1ex}

\noindent Even classes correspond to the black regions of the
black friezes and odd classes to the white regions; the odd class
containing 1 corresponds to the leftmost white region. The
ordinals of (1) register the number of white holes in the black
regions, and those of (2) the number of black disks and the number
of white holes in them.

Composition in $\mbox{\it Frobse}'$ would be defined by reductions
based on the equations of ${\cal L}_\omega'$, like those we gave
for \emph{Frobse}. Essentially, we would have to change only the
reductions corresponding to ${(ab\:3.1)}$, ${(ab\:3.2)}$ and
${(ab\:3.3)}$. We could have instead
\begin{center}
\begin{picture}(310,60)

\put(-13,30){\makebox(0,0)[l]{$(ab\:3.1)$}}
\put(32,31){\makebox(0,0)[l]{$a^n_{2m+1} b^0_{2m+2}=c^n_{2m+2}$}}

\put(157,54){\makebox(0,0)[b]{\tiny $2m\pl 1$}}
\put(180,54){\makebox(0,0)[b]{\tiny $2m\pl 2$}}

\put(167,5){\makebox(0,0)[b]{\tiny $2m\pl 1$}}
\put(190,5){\makebox(0,0)[b]{\tiny $2m\pl 2$}}

\put(170,10){\line(0,1){10}} \put(170,20){\line(-1,1){10}}
\put(170,20){\line(1,1){10}} \put(190,10){\line(0,1){20}}
\put(160,32){\line(0,1){20}} \put(180,52){\line(0,-1){10}}
\put(180,42){\line(-1,-1){10}} \put(180,42){\line(1,-1){10}}

\put(170,29){\circle*{2}} \put(180,33){\circle*{2}}
\put(170,27){\makebox(0,0)[t]{\tiny $n$}}
\put(180,35){\makebox(0,0)[b]{\tiny $0$}}

\put(230,31){\makebox(0,0){$\leadsto$}}

\put(290,54){\makebox(0,0)[b]{\tiny $2m\pl 2$}}

\put(290,5){\makebox(0,0)[b]{\tiny $2m\pl 2$}}

\put(290,10){\line(0,1){42}}

\put(291,32){\makebox(0,0)[l]{\tiny $n$}}

\end{picture}

\begin{picture}(310,60)

\put(32,31){\makebox(0,0)[l]{$a^0_{2m+2} b^n_{2m+3}=c^n_{2m+2}$}}

\put(157,54){\makebox(0,0)[b]{\tiny $2m\pl 2$}}
\put(180,54){\makebox(0,0)[b]{\tiny $2m\pl 3$}}

\put(167,5){\makebox(0,0)[b]{\tiny $2m\pl 2$}}
\put(190,5){\makebox(0,0)[b]{\tiny $2m\pl 3$}}

\put(170,10){\line(0,1){10}} \put(170,20){\line(-1,1){10}}
\put(170,20){\line(1,1){10}} \put(190,10){\line(0,1){20}}
\put(160,32){\line(0,1){20}} \put(180,52){\line(0,-1){10}}
\put(180,42){\line(-1,-1){10}} \put(180,42){\line(1,-1){10}}

\put(170,29){\circle*{2}} \put(180,33){\circle*{2}}
\put(170,27){\makebox(0,0)[t]{\tiny $0$}}
\put(180,35){\makebox(0,0)[b]{\tiny $n$}}

\put(230,31){\makebox(0,0){$\leadsto$}}

\put(290,54){\makebox(0,0)[b]{\tiny $2m\pl 2$}}

\put(290,5){\makebox(0,0)[b]{\tiny $2m\pl 2$}}

\put(290,10){\line(0,1){42}}

\put(291,32){\makebox(0,0)[l]{\tiny $n$}}

\end{picture}

\begin{picture}(310,60)

\put(-13,30){\makebox(0,0)[l]{$(ab\:3.3)$}}
\put(32,31){\makebox(0,0)[l]{$a^n_{2m+1}
b^l_{2m+1}=c^{\omega^{n+l}}_1$}}

\put(177,54){\makebox(0,0)[b]{\tiny $2m\pl 1$}}

\put(177,5){\makebox(0,0)[b]{\tiny $2m\pl 1$}}

\put(180,10){\line(0,1){10}} \put(180,20){\line(-1,1){10}}
\put(180,20){\line(1,1){10}} \put(180,52){\line(0,-1){10}}
\put(180,42){\line(-1,-1){10}} \put(180,42){\line(1,-1){10}}

\put(180,33){\circle*{2}} \put(180,29){\circle*{2}}
\put(180,27){\makebox(0,0)[t]{\tiny $n$}}
\put(180,35){\makebox(0,0)[b]{\tiny $l$}}

\put(230,31){\makebox(0,0){$\leadsto$}}

\put(290,54){\makebox(0,0)[b]{\tiny $1$}}

\put(290,5){\makebox(0,0)[b]{\tiny $1$}}

\put(290,10){\line(0,1){42}}

\put(291,31){\makebox(0,0)[l]{\tiny $\omega^{n+l}$}}

\end{picture}

\begin{picture}(310,60)

\put(32,31){\makebox(0,0)[l]{$a^0_{2m+2} b^0_{2m+2}=c^1_{2m+2}$}}

\put(177,54){\makebox(0,0)[b]{\tiny $2m\pl 2$}}

\put(177,5){\makebox(0,0)[b]{\tiny $2m\pl 2$}}

\put(180,10){\line(0,1){10}} \put(180,20){\line(-1,1){10}}
\put(180,20){\line(1,1){10}} \put(180,52){\line(0,-1){10}}
\put(180,42){\line(-1,-1){10}} \put(180,42){\line(1,-1){10}}

\put(180,33){\circle*{2}} \put(180,29){\circle*{2}}

\put(230,31){\makebox(0,0){$\leadsto$}}

\put(290,54){\makebox(0,0)[b]{\tiny $2m\pl 2$}}

\put(290,5){\makebox(0,0)[b]{\tiny $2m\pl 2$}}

\put(290,10){\line(0,1){42}}

\put(291,31){\makebox(0,0)[l]{\tiny $1$}}

\end{picture}
\end{center}
and analogous reductions for ${(ab\:3.2)}$.

\section{Separable matrix Frobenius monads}

In the preceding section, we saw how symmetry in the category
\emph{Mat} induces a collapse of the ordinals in $\varepsilon_0$
of \emph{Frob} into the ordinals in $\omega^\omega$. In all that,
we have not considered commutative Frobenius monads, which play a
central role in connection with topological quantum field
theories. (For the notion of commutative Frobenius monad, in which
we have a natural transformation from $MM$ to $MM$ with the two
$M$'s ``permuted'', and with appropriate coherence equations, see
\cite{DP09}, Section~3; this notion should not be confused with
the commutative monads of \cite{Ko70}.) With commutative Frobenius
monads, our ordinals are still contained in $\omega^\omega$, as in
the preceding section.

Another collapse of ordinals comes with separability (see
\cite{DI71}, \cite{C91} and \cite{RSW05}). The \emph{separability}
equation for Frobenius monads is the equation
\[
\delta^\Diamond_A\cirk\delta^\Box_A=\mj_{MA}.
\]
If we consider extending \emph{Frob} with this equation, we just
replace $A$ by $n$. To state the consequence of the corresponding
equation ${c^1_{2n+2}=1}$ for ${\cal L}_\omega$, we need some
terminology.

Let the ordinal $0$ be of \emph{even height}. If
$\alpha_1,\ldots,\alpha_n$ are all of \emph{even (odd) height},
then $\omega^{\alpha_1}\sharp\ldots\sharp\,\omega^{\alpha_n}$ is
of \emph{odd (even) height}. If an ordinal in $\varepsilon_0$ is
of even or odd height, we say that it has a \emph{homogeneous
height}. Not all ordinals in $\varepsilon_0$ have a homogeneous
height. The consequence of the separability equation for ${\cal
L}_\omega$ is that every $c^\alpha_n$ is equal to $c^{\alpha'}_n$
for $\alpha'$ an ordinal in $\varepsilon_0$ of homogeneous height;
if $n$ is ${2m\pl 2}$, then $\alpha'$ is of even height, and if
$n$ is ${2m\pl 1}$, then $\alpha'$ is of odd height.

If we combine the separability equation with the equation $(\Phi)$
of the preceding section, then the ordinals in $\varepsilon_0$
collapse to the ordinals in $\omega$. More precisely, the
consequence for ${\cal L}_\omega$ is that we could take as
primitive only the terms $e^k_n$, for $e$ being $a$, $b$ or $c$,
and ${k\in\omega}$, where only $c^k_1$ may have ${k\geq 0}$; in
all other cases, ${k=0}$. In the presence of the separability
equation, the equation
\begin{tabbing}
\hspace{1.7em}$(\Phi^0)$\hspace{12em}$\Phi^0_{MA}=M\Phi^0_A$
\end{tabbing}
has the same force as the equations $(\Phi)$. According to our
definition, $\Phi^0_A$ is
${\varepsilon^\Box_A\cirk\varepsilon^\Diamond_A}$.

We call Frobenius monads that satisfy $(\Phi^0)$ and the
separability equation \emph{separable matrix} Frobenius monads.
For separable matrix Frobenius monads, we can answer positively
the question of sufficiency left open at the beginning of the
preceding section. Namely, there is a faithful monoidal functor
$F$ from the separable matrix Frobenius monad generated by a
single object into the category \emph{Mat}. In fact, something
stronger holds: for every natural number ${p\geq 2}$, there is a
functor $F$ as above such that ${F(1)=p}$. We will not prove this
in detail, but just give some indications.

Our task is to represent in \emph{Mat} an ordered pair made of a
maximal split equivalence (see Section~6) and a natural number,
which is the ordinal ${k\in\omega}$ tied to $c^k_1$. We may reject
the odd equivalence classes from this maximal split equivalence,
and then represent the remaining split equivalence in a Brauerian
manner (see \cite{DP03}, \cite{DP03a} and \cite{DP06}). The
natural number $k$ will be mapped to the scalar $p^k$. This is
analogous to representing ${\cal K}_c$ in \emph{Mat} (in the
section with that name in \cite{DP03}), but is not exactly the
same. In the free self-adjunction ${\cal K}_c$ of the $\cal K$
type (corresponding to Temperley-Lieb algebras), the ordinals in
$\varepsilon_0$ of ${\cal L}_\omega$ also collapse to natural
numbers, and are not tied to particular regions of the frieze.
This is analogous to what we have with separable matrix Frobenius
monads, but is not exactly the same. The difference is that for
${\cal K}_c$ all circles are counted, while here we count circles
tied to ${\varepsilon^\Box_A\cirk\varepsilon^\Diamond_A}$, which
may be moved according to the equation $(\Phi)$ or $(\Phi^0)$, and
do not count circles tied to
${\delta^\Diamond_A\cirk\delta^\Box_A}$, according to the
separability equation. We will deal with these matters in more
detail on another occasion.

Let us sum up matters from the preceding section and the present
one. We know that the equation $(\Phi)$ is necessary for the
existence of a faithful monoidal functor $F$ into the category
\emph{Mat}. We do not know whether $(\Phi)$ is sufficient. If it
were, then we could legitimately call Frobenius monads that
satisfy $(\Phi)$ \emph{matrix} Frobenius monads. We know on the
other hand that $(\Phi)$ together with the separability equation
is sufficient for the existence of such an $F$, but we do not know
whether the separability equation is necessary, though this
necessity does not seem likely. Since ordinals in separable matrix
Frobenius monads have collapsed to natural numbers, with these
monads we reach the boundary we set ourselves for this paper,
where we wanted to investigate the role of bigger ordinals in
Frobenius monads.

\vspace{2ex}

\noindent {\small {\it Acknowledgement$\,$}. Work on this paper
was supported by the Ministry of Science of Serbia (Grant ON174026).}

\end{document}